\documentclass[12pt]{article}

\usepackage{euscript}

\usepackage{amsfonts}
\usepackage{amsmath}
\usepackage{amssymb}
\usepackage{latexsym,color}
\usepackage{amsbsy}
\usepackage{graphicx}

\newtheorem{thm}{Theorem}[section]
\newtheorem{con}{Conjecture}[section]
\newtheorem{lemma}[thm]{Lemma}
\newtheorem{cor}[thm]{Corollary}
\newtheorem{pro}[thm]{Proposition}
\newtheorem{example}[thm]{Example}
\newtheorem{definition}[thm]{Definition}
\newcommand{\tab}{\mbox{tab}}
\newtheorem{remark}[thm]{Remark}
\newtheorem{Algorithm}[thm]{Algorithm}

\newcommand{\pc}{{\cal P}}

\newcommand{\comment}[1]{}

\newcommand{\Bx}{B_X}

\newcommand{\C}{{\mathbb C}}   
\newcommand{\N}{{\mathbb N}}   
\newcommand{\Po}{{\mathbb P}}   
\newcommand{\Zi}{{\mathbb Z}}

\newcommand{\da}{\downarrow}
\newcommand{\spec}{\mbox{spec}_G}
\newcommand{\cont}{\mbox{cont}}
\newcommand{\contg}{\mbox{cont}_G}
\newcommand{\tabx}{\mbox{tab}_G}

\newcommand{\I}{{\cal{I}}}
\newcommand{\J}{{\cal{J}}}
\newcommand{\R}{{\cal{R}}}

\newcommand{\ncom}{\newcommand}
\ncom{\ns}{\normalsize}
\ncom{\la}{\lambda}
\ncom{\bm}{\boldmath}
\ncom{\noi}{\noindent}
\ncom{\bq}{\begin{equation}}
\ncom{\eq}{\end{equation}}  
\ncom{\beqn}{\begin{eqnarray*}}
\ncom{\eeqn}{\end{eqnarray*}}  
\ncom{\ba}{\begin{array}}

\ncom{\ea}{\end{array}}
\ncom{\beq}{\begin{eqnarray}}
\ncom{\eeq}{\end{eqnarray}}
\ncom{\nno}{\nonumber}
\ncom{\hs}{\mbox{\hspace{.25cm}}}
\ncom{\rar}{\rightarrow}
\ncom{\Rar}{\Rightarrow}  
\ncom{\noin}{\noindent}   
\ncom{\bc}{\begin{center}}
\ncom{\ec}{\end{center}}  
\ncom{\sz}{\scriptsize}   
\ncom{\fpd}{\Phi(\pi^{'})}
\ncom{\fp}{\Phi(\pi) }
\ncom{\nk}{\left< \begin{array}{c}
                       n\\k \end{array} \right>}
\ncom{\nd}{1^{'},2^{'},\cdots,n^{'}}
\ncom{\de}{\bigtriangleup (F_{2n},\leq)}
\ncom{\del}{\bigtriangleup} 
\ncom{\cov}{<\!\!\!\!\cdot }

\ncom{\bt}{\begin{thm}}
\ncom{\bcon}{\begin{con}}
\ncom{\et}{\end{thm}}
\ncom{\econ}{\end{con}}
\ncom{\bl}{\begin{lemma}}
\ncom{\el}{\end{lemma}}  
\ncom{\bco}{\begin{cor}} 
\ncom{\ds}{\displaystyle}
\ncom{\eco}{\end{cor}}   
\ncom{\bp}{\begin{pro}}  
\ncom{\ep}{\end{pro}}    
\ncom{\bex}{\begin{example}}
\ncom{\seq}{\subseteq}
\ncom{\eex}{\end{example}}  
\ncom{\bd}{\begin{definition}}
\ncom{\ed}{\end{definition}}  
\ncom{\brm}{\begin{remark}}   
\ncom{\erm}{\end{remark}}     
\ncom{\bal}{\begin{Algorithm}}
\ncom{\eal}{\end{Algorithm}}  
\ncom{\ol}{\overline}

\ncom{\pf}{\noi {\bf Proof  }}
\ncom{\be}{\begin{enumerate}} 
\ncom{\ee}{\end{enumerate}}   
\ncom{\s}{\subset}
\ncom{\T}{{\cal T}}
\ncom{\B}{{\cal B}}
\ncom{\A}{{\cal A}}
\ncom{\Y}{{\cal Y}}

\title{\Large{{\textcolor{black} {\bf On the representation theory of $G\sim
S_n$}}}}

\author{{\textcolor{black} { Ashish Mishra and Murali K. Srinivasan}} \\
{\em  {Department of Mathematics}}\\
{\em  {Indian Institute of Technology, Bombay}}\\
{\em  {Powai, Mumbai 400076, INDIA}}\\
{\bf  \texttt{ashishm@math.iitb.ac.in}}\\ 
{\bf  \texttt{murali.k.srinivasan@gmail.com}}}

\begin{document}
\date{}
\maketitle

\begin{abstract} In the Vershik-Okounkov approach to the
complex irreducible representations of $S_n$ and $G\sim S_n$ we
parametrize the irreducible representations and their bases 
by spectral objects rather than
combinatorial objects and then, at the end, give a bijection between the
spectral and combinatorial objects. 
The fundamental ideas are similar in both cases but there are
additional technicalities involved in the $G\sim S_n$ case.
This was carried out by Pushkarev. 

The present work gives a fully 
detailed exposition of Pushkarev's theory. For the most part we follow the
original but our definition of a Gelfand-Tsetlin subspace,
based on a multiplicity free chain of subgroups, 
is slightly different and leads to a more natural development of the theory.
We also work out in
detail an example, the generalized Johnson scheme, from 
this viewpoint. 
\end{abstract}

\section{\large{Introduction}} 

Let $G$ be a finite group. The symmetric group $S_n$ acts on $G^n=G\times
\cdots \times G$ ($n$ factors) by permuting the coordinates and this action
defines the semidirect product $G^n \rtimes S_n$ of $G^n$ by $S_n$.
The group $G^n\rtimes S_n$ is called the {\em wreath product} of $G$ by
$S_n$ and is denoted $G\sim S_n$ (our notation follows {\bf\cite{m}}).
We set $G_n=G\sim S_n$. The elements of $G_n$ are the set of all
$(n+1)$-tuples $(g_1,\ldots ,g_n,\pi)$, where $\pi\in S_n$, and $g_i\in G$
for all $i$. The multiplication rule and inverse of an element in $G_n$ are
given by
\beqn
(g_1,\ldots ,g_n,\pi)(h_1,\ldots ,h_n,\tau)&=&(g_1h_{\pi^{-1}(1)},\ldots
,g_nh_{\pi^{-1}(n)}, \pi\tau),\\
(g_1,\ldots ,g_n,\pi)^{-1}&=& (g_{\pi(1)}^{-1},\ldots
,g_{\pi(n)}^{-1},\pi^{-1}).
\eeqn

The complex representation theory of $G_n$ is a
classical and well studied topic. Among the many sources we mention James
and Kerber {\bf\cite{jk}}, Macdonald {\bf\cite{m}}, and the recent book of
Ceccherini-Silberstein, Scarabotti, and Tolli {\bf\cite{cst3}}. The basic
problem can be stated as follows.

Let $\pc$ denote the set of all partitions (there is a unique partition of
zero with zero parts) and let $\pc_n$ denote the set of all partitions of
$n$. For a finite set $X$, we define 
$$\pc(X) =  \{ \mu \;|\; \mu: X \rar \pc\}.$$
For $\mu\in \pc(X)$, define $\| \mu \| = \sum_{x\in X} |\mu(x)|$, where
$|\mu(x)|$ is the sum of the parts of the partition $\mu(x)$ and define
$$\pc_n(X)=\{\mu\in \pc(X) \;|\; \| \mu \|=n\}.$$  

Let $G_*$ denote the set of conjugacy classes in $G$. The conjugacy classes
of $G_n$ are parametrized by $\pc_n(G_*)$ ({\bf\cite{jk,m,cst3}}).

Let $\Y$ denote the set of all Young diagrams (there is a unique Young
diagram with zero boxes) and $\Y_n$ denote the set of all Young diagrams
with $n$ boxes. For a finite set $X$, we define 
$$\Y(X) =  \{ \mu \;|\; \mu: X \rar \Y\}.$$
For $\mu\in \Y(X)$, define $\| \mu \| = \sum_{x\in X} |\mu(x)|$, where
$|\mu(x)|$ is the number of boxes of the Young diagram $\mu(x)$ and define
$$\Y_n(X)=\{\mu\in \Y(X) \;|\; \| \mu \|=n\}.$$

Denote by $G^{\wedge}$ the (finite) set of equivalence classes of finite
dimensional complex irreducible representations of $G$. Given $\sigma \in
G^{\wedge}$, we denote by $V^{\sigma}$ the corresponding irreducible
$G$-module. Elements of
$\Y(G^{\wedge})$ are called {\em Young $G$-diagrams} and elements of
$\Y_n(G^{\wedge})$ are called {\em Young $G$-diagrams with $n$ boxes}. 
Given $\mu\in \Y(G^{\wedge})$ and $\sigma \in G^{\wedge}$, we denote by 
$\mu \da \sigma$ the set of all Young $G$-diagrams obtained from $\mu$ by
removing one of the inner corners in the Young diagram $\mu(\sigma)$.

Let $\mu\in \Y$. A {\em Young tableau of shape $\mu$} 
is obtained by taking the Young diagram $\mu$ and filling its 
$|\mu|$ boxes 
(bijectively) with the numbers $1,2,\ldots ,| \mu |$. 
A Young tableau  is said to be {\em standard} if the numbers in
the boxes strictly increase along each row and each column of the Young
diagram of $\mu$. Let $\tab(n,\mu)$, where $\mu \in
\Y_n$,
denote the set of all
standard Young tableaux of shape $\mu$ and let $\tab(n)=\cup_{\mu\in
\Y_n}\tab(n,\mu)$.  

Let $\mu\in \Y(G^{\wedge})$. A {\em Young $G$-tableau of shape $\mu$} 
is obtained by taking the Young $G$-diagram $\mu$ and filling its 
$\|\mu\|$ boxes 
(bijectively) with the numbers $1,2,\ldots ,\| \mu \|$. 
A Young $G$-tableau   is said to be {\em standard} if the numbers in
the boxes strictly increase along each row and each column of all Young
diagrams occuring in $\mu$. Let $\tabx(n,\mu)$, where $\mu \in
\Y_n(G^{\wedge})$,
denote the set of all
standard Young $G$-tableaux of shape $\mu$ and let $\tabx(n)=\cup_{\mu\in
\Y_n(G^{\wedge})}\tabx(n,\mu)$.  

Let $T\in \tabx(n)$ and $i\in \{1,\ldots ,n\}$.  
If $i$ 
appears in the Young diagram
$\mu(\sigma)$, where $\mu$ is the shape of $T$ and 
$\sigma\in G^{\wedge}$, we write $r_T(i)=\sigma$.

The complex 
irreducible representations of $G_n$ are parametrized by $\Y_n(G^{\wedge})$ 
and the basic problem of the representation theory of $G_n$ is to
explain this correspondence between irreducible representations of $G_n$
and elements
of $\Y_n(G^{\wedge})$. This is done in {\bf\cite{m}} using symmetric
functions and the characteristic map and in {\bf\cite{jk,cst3}} using
Clifford theory and the little group method.

In {\bf\cite{p}} Pushkarev, building on the Vershik-Okounkov approach in the
$S_n$ case {\bf\cite{vo1,vo2,cst2}}, gave a spectral explanation for this
correspondence, namely, an internal analysis of the irreducible
representations of $G_n$ yields spectral objects parametrizing the
irreducible representations and then
a bijection is given between these spectral objects and $\Y_n(G^{\wedge})$. 
This approach is inductive in nature and has the following advantages:

\noi (a) The group $G_n$ can be identified with the
subgroup 
$$\{ (g_1,\ldots ,g_n,e,\pi) \;|\; \pi\in S_{n+1}\mbox{ with }\pi(n+1)=n+1
\mbox{ and }g_i\in G,\;1\leq i \leq n\}$$
of $G_{n+1}$ ($e = $ identity element of $G$) and we have an infinite chain of finite groups
$$G_1\seq G_2 \seq \cdots .$$
As a natural byproduct of the theory we get the branching rule from
$G_{n+1}$ to $G_n$: denote the irreducible $G_{n+1}$-module corresponding to
$\mu\in \Y_{n+1}(G^{\wedge})$ by $V^{\mu}$. Then we have $G_n$-module
isomorphisms
\beqn
V^{\mu} &\cong& \oplus_{\sigma\in G^{\wedge}} \dim(V^{\sigma})\;\left( 
\oplus_{\lambda \in \mu \da \sigma} V^{\lambda} \right).
\eeqn

\noi (b) Another natural byproduct of the theory yields a parametrization of
the bases of irreducible $G_n$-modules using standard Young $G$-tableaux 
and bases of irreducible $G^n$-modules. More precisely, for $\mu\in
\Y_n(G^{\wedge})$, we have a canonical direct sum decomposition of $V^{\mu}$ 
into
subspaces, called Gelfand-Tsetlin subspaces, 
\beqn V^{\mu} &=& \oplus_{T\in \tabx(n,\mu)} V_T,
\eeqn
where each $V_T$ is closed under the action of $G^n=G\times \cdots \times G$
($n$ factors) and, as a $G^n$-module, is isomorphic to the irreducible
$G^n$-module
$$V^{r_T(1)}\otimes V^{r_T(2)}\otimes \cdots \otimes V^{r_T(n)}.$$  

The present work gives a fully detailed exposition
of Pushkarev's theory. Our development, based on a multiplicity free chain
of subgroups, is slightly
different from the original and is along the following lines.

For
$g\in G$ and $1\leq i \leq n$ we denote by $g^{(i)}$ the element 
$(e,\ldots ,e,g,e,\ldots ,e,1)\in G_n$, where $g$ is in the $i$th
spot, $e$ denotes the identity element of $G$, and $1$ denotes the identity
element of $S_n$. Denote by $G^{(i)}$ the subgroup $\{g^{(i)} | g\in G \}$  
of $G_n$. Note that $G^{(1)},\ldots ,G^{(n)}$ commute. 
We may also think of $S_n$ as the subgroup $\{(e,\ldots ,e,\pi) | \pi\in
S_n\}$. We write the element $(e,\ldots ,e,\pi)$ as $\pi$.
We may thus write an
element $(g_1,\ldots ,g_n,\pi)\in G_n$ as
$g_1^{(1)}\ldots g_n^{(n)}\pi = 
\pi g_1^{(\pi^{-1}(1))}\ldots
g_n^{(\pi^{-1}(n))}= 
\pi g_{\pi(1)}^{(1)}\ldots
g_{\pi(n)}^{(n)}$. 

For $n\geq 1$, set
$H_{n,n}=G_n $ and consider the following chain of subgroups
\beq \label{mfic}
&H_{1,n}\seq H_{2,n} \seq \cdots \seq H_{n,n},&
\eeq
where, for $1\leq i \leq n$,
$$ H_{i,n} = \{(g_1,\ldots ,g_n,\pi)\in G_n \;|\;
            \pi(j)=j \mbox{ for }i+1\leq j\leq n\}.  
$$
Note that $H_{1,n}$ is isomorphic to $G^n$. The following are the main steps
in the representation theory of $G_n$. 

\noi (i) A direct argument shows that branching from
$H_{i,n}$ to $H_{i-1,n}$ is simple, i.e., multiplicity free.

\noi (ii) Consider an irreducible $H_{m,n}$-module $V$.
Since the branching is
simple the decomposition of $V$ into irreducible $H_{m-1,n}$-modules is
canonical. Each of these modules, in turn, decompose canonically into  
irreducible $H_{m-2,n}$-modules. Iterating this construction we get a  
canonical decomposition of $V$ into irreducible $G^n=H_{1,n}$-modules,
called
the {\em Gelfand-Tsetlin decomposition (GZ-decomposition)} of $V$. 
The irreducible $G^n$-modules in this decomposition are called the 
{\em Gelfand-Tsetlin subspaces (GZ-subspaces)} of $V$.  

\noi (iii) Let $Z_{m,n}$ denote the center of the group algebra $\C[H_{m,n}]$.
The {\em Gelfand-Tsetlin algebra (GZ-algebra)}, denoted $GZ_{m,n}$,
is defined to be
the (commutative) subalgebra of $\C[H_{m,n}]$
generated by $Z_{1,n}\cup Z_{2,n}\cup \cdots  \cup Z_{m,n}$.
It is shown that $GZ_{m,n}$ consists of
all elements in $\C[H_{m,n}]$ that act by scalars on the
GZ-subspaces in every irreducible representation of $H_{m,n}$.
It follows that if we have a finite generating set for $GZ_{m,n}$ then the
GZ-subspaces are determined by the eigenvalues on this generating set.

\noi (iv) Following Pushkarev, for $i=1,2,\ldots ,n$,
we define the (generalized) YJM elements 
$X_1,X_2,\ldots ,X_n$ of $\C[H_{n,n}]$:  
\beqn
X_i &=& \sum_{k=1}^{i-1}\sum_{g\in G}(g^{-1})^{(k)}g^{(i)} (k,i).
\eeqn
Note that $X_1=0$. For an algebra $A$, let $Z[A]$ denote the center of $A$.
It is shown that
$GZ_{m,n} = \langle Z[\C[G^n]], X_1,X_2,\ldots ,X_m \rangle $.

\noi (v) 
By a GZ-subspace of $G_n$ we mean a GZ-subspace in some
irreducible representation of $G_n$. Let $W$ be a GZ-subspace of
$G_n$. Then $W$ is an irreducible $G^n$-module and hence is isomorphic
to $V^{\rho_1}\otimes \cdots \otimes V^{\rho_n}$, where $\rho_i \in
G^{\wedge}$, for all $i$. We call $\rho = (\rho_1,\ldots ,\rho_n)$ the {\em
label} of the GZ-subspace $W$.

It follows from steps (iii) and  (iv) above 
that a GZ-subspace $W$ of $G_n$ is uniquely determined by
its label and the eigenvalues of $X_1,\ldots ,X_n$ on $W$.
To a GZ-subspace $W$ we associate the tuple 
$$\alpha(W) =
(\rho,a_1,a_2,\ldots ,a_n),$$
where $\rho$ is the label of $W$ and
$a_i = \mbox{ eigenvalue of $X_i$  on } W$. We call $\alpha(W)$ the {\em
weight} of the GZ-subspace $W$. 
Define
$$\spec (n) = \{\alpha (W)\;:\;W \mbox{ is a GZ-subspace of $G_n$}\},$$
called the spectrum of $G_n$.

We have $\mbox{dim }GZ_{n,n} = |\spec(n)|$. There is a natural equivalence
relation $\sim$ on $\spec(n)$: for $\alpha, \beta \in \spec(n)$, $\alpha  
\sim \beta$ iff the corresponding GZ-subspaces are in same
$G_n$-irreducible. Clearly, we have $|\spec(n)/\sim| = |G_n^{\wedge}|$.

The representation theory of $G_n$ is governed by the spectral object
$\spec(n)$.

\noi (vi) In the final step we construct a bijection between $\spec(n)$ and
$\tabx(n)$ such that tuples in $\spec (n)$ 
related by $\sim$ go to standard Young $G$-tableaux of the
same shape. This step is carried out inductively using an 
analysis of the following commutation relations
that hold in $G_n$ (where $s_i=\mbox{ the Coxeter generator }(i,
i+1)$):

(a) $X_1,\ldots ,X_n$ commute.

(b) $X_ig^{(l)}=g^{(l)}X_i$, $g\in G, 1\leq i,l \leq n$.

(c) $s_ig^{(i)}s_i = g^{(i+1)}$, $g\in G, 1\leq i \leq n-1$. In
particular,
$s_i^2=1$. 

(d) $s_ig^{(l)}=g^{(l)}s_i$, $1\leq i\leq n-1, 1\leq l \leq n, l\not=
i,i+1$.

(e) $s_iX_is_i + \sum_{g\in G} g^{(i+1)}s_i(g^{-1})^{(i+1)}=X_{i+1}$,
$1\leq i \leq n-1$.

(f)  $s_iX_l=X_ls_i$,  $1\leq i\leq n-1, 1\leq l \leq n, l\not= i,i+1$.

We now give a brief synopsis of the paper. Section 2 collects some
preliminaries on wreath products. Section 3 discusses Gelfand-Tsetlin
subspaces, Gelfand-Tsetlin decompositions, and Gelfand-Tsetlin algebras for
an inductive chain of finite groups with simple branching. In Section 4 we
first show that the chain (\ref{mfic}) is multiplicity free and then show
that the corresponding Gelfand-Tsetlin algebras are generated over
$Z[\C[G^n]]$ by the YJM elements, thereby defining the weight of a
GZ-subspace and the spectrum $\spec(n)$ of $G_n$. 
Section 5 describes, using the commutation relations (a)-(f) above, 
the action of the
Coxeter generators on the Gelfand-Tsetlin subspaces in terms of
transformations of weights. In Section 6, using the results of Section 5, 
we give a bijection between $\spec(n)$ and $\tabx(n)$ via the content
vectors of standard Young $G$-tableaux. In Section 7 we study the simplest
nontrivial example of the Vershik-Okounkov theory, 
the classical ``Johnson schemes" and the ``generalized Johnson schemes". 
We consider multiplicity
free $S_n$, $G_n$-actions and explicitly write down the GZ-vectors (in the
$S_n$ case) and the GZ-subspaces (in the $G_n$ case) and also identify the
irreducibles which occur. 

\section{\large{Preliminaries}} 

The positive integers are denoted $\Po$ and the nonnegative integers are
denoted $\N$.

We enumerate the conjugacy classes of $G$ as $G_*=\{C_1,\ldots ,C_t\}$ and
assume that $C_1=\{e\}$. 
We say that $g\in G$ is of {\em type} $j$ if $g\in C_j$.
Define an involution $\I :\{1,\ldots ,t\} \rar
\{1,\ldots ,t\}$ as follows: $\I(j)=j'$ if $j'$ is the type of $g^{-1}$, for
$g\in C_j$.

Let $h=(g_1,\ldots ,g_n,\pi)\in G_n$ and let $\tau = (i_1,i_2,\ldots
,i_k)$ be a $k$-cycle in $\pi$. The element $g_{i_k}g_{i_{k-1}}\cdots
g_{i_1}\in G$ is called the {\em cycle product} of 
$h$ corresponding to the cycle
$\tau$ of $\pi$ and its type is easily seen to be independent of the order
in which the elements of $\tau$ are listed. 
Thus we may define $\rho_h:G_* \rar \pc$ by
$$\rho_h(C_i) = \mbox{ Multiset of lengths of all cycles of $\pi$ whose 
cycle product lies in }C_i,\;1\leq i \leq t.  
$$
Clearly $\rho_h\in \pc_n(G_*)$. We say that $\rho_h$ is the {\em type} of
$h\in G_n$. 

Suppose two elements $(g_1,\ldots ,g_n,\pi)$ and $(f_1,\ldots ,f_n,\tau)$
are conjugate in $G_n$. Then we have
\beq \label{co1}
(h_1,\ldots ,h_n,\sigma)(g_1,\ldots ,g_n,\pi)(h^{-1}_{\sigma(1)},\ldots
,h^{-1}_{\sigma(n)},\sigma^{-1}) &=& (f_1,\ldots ,f_n,\tau),
\eeq
for some $h_1,\ldots ,h_n\in G$ and $\sigma \in S_n$. Thus $\tau = \sigma
\pi \sigma^{-1}$ and $\tau$ and $\pi$ are conjugate in $S_n$. 

We now want to consider the cycle products in $(g_1,\ldots ,g_n,\pi)$ and
$(f_1,\ldots ,f_n,\tau)$. For simplicity we shall write the element
$(g_1,\ldots ,g_n,\pi)$ as $(\ldots,g_i,\ldots,\pi)$ (it being understood
that $g_i$ is in the $i$th spot). We have 
\beqn
\lefteqn{(\ldots,h_i,\ldots,\sigma)(\ldots,g_i,\ldots,\pi)
(\ldots,h^{-1}_{\sigma(i)},\ldots,\sigma^{-1})}\\
&=&
(\ldots,h_ig_{\sigma^{-1}(i)},\ldots,\sigma\pi)
(\ldots,h^{-1}_{\sigma(i)},\ldots,\sigma^{-1})\\
&=&
(\ldots,h_ig_{\sigma^{-1}(i)}
h^{-1}_{\sigma\pi^{-1}\sigma^{-1}(i)},\ldots,\sigma\pi\sigma^{-1})\\
&=&
(\ldots,h_ig_{\sigma^{-1}(i)}
h^{-1}_{\tau^{-1}(i)},\ldots,\tau)\\
&=&
(\ldots,f_i,\ldots,\tau)\\
\eeqn
Let $(i_1,\ldots ,i_k)$ be a cycle in $\pi$. Then $(\sigma(i_1),\ldots
,\sigma(i_k))$ is a cycle in $\tau$. We have, using the calculation above,
\beq
f_{\sigma(i_k)}&=& h_{\sigma(i_k)}g_{\sigma^{-1}(\sigma(i_k))}
                   h^{-1}_{\tau^{-1}(\sigma(i_k))}\\ \label{co2}
               &=& h_{\sigma(i_k)}g_{i_k}h^{-1}_{\sigma(i_{k-1})}.
\eeq
Thus we have (using $\tau^{-1}(\sigma(i_1))=\sigma(i_k)$)
\beqn
\lefteqn{f_{\sigma(i_k)}f_{\sigma(i_{k-1})}\cdots f_{\sigma(i_1)} }\\
&=&(h_{\sigma(i_k)}g_{i_k}h^{-1}_{\sigma(i_{k-1})})
(h_{\sigma(i_{k-1})}g_{i_{k-1}}h^{-1}_{\sigma(i_{k-2})})\cdots
(h_{\sigma(i_{1})}g_{i_{1}}h^{-1}_{\sigma(i_{k})})\\
&=&h_{\sigma(i_k)}g_{i_k}\cdots g_{i_1}h^{-1}_{\sigma(i_k)}.
\eeqn
Thus the type of the cycle products $g_{i_k}\cdots g_{i_1}$ and
$f_{\sigma(i_k)}\cdots f_{\sigma(i_1)}$ are the same. It follows that if
two elements of $G_n$ are conjugate then they  have the same type. 

Conversely, suppose that $(g_1,\ldots ,g_n,\pi), (f_1,\ldots ,f_n,\tau)\in
G_n$ have the same type. Then we can easily write down a $\sigma \in S_n$
such that $\sigma \pi \sigma^{-1}=\tau$ and such that, for every cycle
$(i_1,\ldots ,i_k)$ of $\pi$, the cycle products $g_{i_k}\cdots g_{i_1}$ and
$f_{\sigma(i_k)}\cdots f_{\sigma(i_1)}$ have the same type. Now, using
(\ref{co2}), we can find $h_1,\ldots ,h_n\in G$ such that (\ref{co1}) holds.
It follows that two elements of $G_n$ are conjugate if and only if they have
the same type.

An element $g=(g_1,\ldots ,g_n,\pi)\in G_n$ is said to be a {\em
nontrivial cycle of type $j$} if (exactly) one of the following
conditions hold:

(i) All cycles of $\pi$ have length 1 (i.e., $\pi$ is the identity
permutation) and, for some $1\leq i\leq n$,
$g_l=e$ for $l\not=i$, $g_i$ is of type $j$, and $2\leq j \leq t$.
We say that $\{i\}$ is the {\em support} of $g$. We say that $g$ is a {\em
nontrivial $1$-cycle} of type $j$.

(ii) There is exactly one cycle, say $(i_1,\ldots ,i_k)$, in the cycle
decomposition of $\pi$ of length $\geq 2$, 
the cycle product $g_{i_k}\ldots g_{i_1}$ is of
type $j$, and $g_l=e$, for $l\not\in \{i_1,\ldots ,i_k\}$. 
We say that $\{i_1,\ldots ,i_k\}$ is the {\em support} of $g$. Note that in
this case there is no restriction on $j$, i.e., $1\leq j \leq t$.
We say that $g$ is a {\em
nontrivial $k$-cycle} of type $j$.

Just as in the $S_n$ case every element of $G_n$ can be written as a
product of commuting nontrivial cycles with disjoint support.

By a {\em nontrivial part} of a partition we mean a part $\geq 2$.
For a partition $\mu\in \pc$ we denote by $\#\mu$ the sum of all the
nontrivial parts (with multiplicity) of $\mu$.

Let $\rho \in \pc_n(G_*)$. By a {\em part} of $\rho$ we mean a pair
$(k,j)$, where $k \in \Po$, $j\in\{1,\ldots ,t\}$, and
$k$ is a part of $\rho(C_j)$. We may specify $\rho$ by giving its multiset of parts
(for example, if $k$ appears $m$ times in $\rho(C_j)$ then the part
$(k,j)$ appears $m$ times in the multiset of parts). 
We say the part $(k,j)$ is {\em nontrivial} if $(k,j)\not= (1,1)$.
We define 
\beqn
\#\rho&=& \sum_{j=2}^t |\rho(C_j)| + \#(\rho(C_1)),
\eeqn
i.e., $\#\rho$ is the sum of the first components (with multiplicity)
of all the nontrivial parts of $\rho$.

For a permutation $s\in S_n$ we denote by $\ell(s)$ the number of inversions
in $s$. It is well known that $s$ can be written as a product of
$\ell(s)$ Coxeter transpositions $s_i = (i, i+1),\;i=1,2,\ldots ,n-1$ and
that $s$ cannot be written as a product of fewer Coxeter transpositions.

All our algebras are finite dimensional, over $\C$, and have units.
Subalgebras contain the unit, and algebra
homomorphisms preserve units.
Given elements or subalgebras $A_1,A_2,\ldots ,A_n$ of an algebra $A$ we
denote by $\langle A_1, A_2, \ldots ,A_n \rangle$ the subalgebra of $A$
generated by $A_1\cup A_2\cup \cdots \cup A_n$.

If $A$ is an algebra and $\rho : A \rar \mbox{End}(V)$ is a representation
then we use several notations for the action of $A$ on the elements of $V$.
For $a\in A$ and $v\in V$ we set
$$\rho(a)(v)=a\cdot v = av = a(v).$$

Similarly, for $a\in A$ and $W\seq V$ we set
$$\rho(a)(W)=a\cdot W = aW = a(W).$$

\section{\large{Gelfand-Tsetlin subspaces, Gelfand-Tsetlin decomposition,
and Gelfand-Tsetlin algebras}} 

The fundamental building blocks of the spectral
approach to the representation theory of $S_n$ and $G_n$ 
are the concepts of Gelfand-Tsetlin subspaces (GZ-subspaces), 
Gelfand-Tsetlin decompositions 
(GZ-decompositions), and
Gelfand-Tsetlin algebras (GZ-algebras), together with a convenient set of
generators for the GZ algebras, for an inductive chain of finite
groups with simple branching. We discuss this in the present and next
sections.

Let
\beq \label{c}
& F_1 \seq F_2 \seq \cdots \seq F_n &
\eeq
be an inductive chain of finite groups. Note that we have not assumed that
$F_1$ is the trivial group with one element. We call $F_1$ the {\em base
group}. 
Define the following directed graph, called the {\em branching
multigraph or Bratelli diagram} of this chain: its vertices are the elements
of the set
$${\ds \coprod_{i = 1}^n} F_i^{\wedge}\;\;\;\mbox{(disjoint union)}$$
and two
vertices $\mu, \lambda$ are joined by $k$ directed edges from $\mu$ to
$\lambda$ whenever $\mu\in F_{i-1}^{\wedge}$ and  $\lambda\in F_i^{\wedge}$
for
some $i$, and the multiplicity of $\mu$ in the restriction of $\lambda$ to
$F_{i-1}$ is $k$. We say that $F_i^{\wedge}$ is {\em level  $i$} of
the branching multigraph.
We write $\mu \nearrow \lambda$ if there is an edge from $\mu$ to $\lambda$.

For the rest of this section assume that the branching multigraph defined
above is actually a graph, i.e., the multiplicities of all restrictions are
0 or 1. We say that the {\em branching or multiplicities} are {\em simple}.

Consider the $F_n$-module $V^{\lambda}$, where $\lambda\in F_n^{\wedge}$.
Since the branching is simple, the decomposition
$$V^{\lambda}=\bigoplus_{\mu} V^{\mu},$$
where the sum is over all $\mu\in F_{n-1}^{\wedge}$ with $\mu\nearrow
\lambda$, is canonical. Iterating this decomposition we obtain a canonical
decomposition of $V^{\lambda}$ into irreducible $F_1$-modules, i.e.,
\beq \label{decomp}
&V^{\lambda}=\bigoplus_{T} V_T,&\eeq
where the sum is over all possible chains
\beq \label{ch}
T&=&\lambda_1\nearrow\lambda_2 \nearrow \cdots \nearrow \lambda_n,
\eeq
with $\lambda_i\in F_i^{\wedge}$ and $\lambda_n = \lambda$.

We call (\ref{decomp}) the {\em Gelfand-Tsetlin
decomposition (GZ-decomposition)} of $V^{\lambda}$ and we call each $V_T$ in
(\ref{decomp}) a {\em Gelfand-Tsetlin subspace (GZ-subspace)} of $V^{\lambda}$. 
By the definition of $V_T$, we have, for $v_T\in V_T$,
\beqn
\C[F_i]\cdot v_T &=& V^{\lambda_{i}},\;\;\;i=1,2,\ldots ,n.
\eeqn
Also note that chains
in (\ref{ch}) are in bijection with directed paths in the branching graph
from an element $\lambda_1$ of $F_1^{\wedge}$ to $\lambda$.

Fix a {\em distinguished basis} $B^{\mu}$ for each $V^{\mu},
\lambda\in F_1^{\wedge}$.
Considering the algebra isomorphism
\beq\label{iso}
\C[F_n]&\cong&\bigoplus_{\lambda\in F_n^{\wedge}} \mbox{End}(V^{\lambda}),
\eeq
given by
$$g \mapsto ( V^{\lambda} \buildrel {g}\over \rightarrow V^{\lambda}\;:\;
\lambda
\in F_n^{\wedge}),\;\;g\in F_n,$$
we can define three natural subalgebras of $\C[F_n]$ based on
the GZ-decomposition (\ref{decomp}).
\beqn
\A_0(n) &=& \{ a\in \C[F_n] : a \mbox{ acts by a scalar on each GZ-subspace of
} V^{\lambda},\;\mbox{for all }\lambda\in F_n^{\wedge}\},\\
\A_1(n) &=& \{ a\in \C[F_n] : a \mbox{ acts diagonally in the distinguished
basis $B^{\lambda}$ of each}\\
&&\;\;\;\;\;\;\;\;\;\;\;\;\;\;\;\;\;\;\;\mbox{GZ-subspace of
} V^{\lambda},\;\mbox{for all }\lambda\in F_n^{\wedge}\},\\
\A_2(n) &=& \{ a\in \C[F_n] :  \mbox{each GZ-subspace of
} V^{\lambda} \mbox{ is $a$ invariant},\;
\mbox{for all }\lambda\in F_n^{\wedge}\}.
\eeqn  
Clearly, $\A_0(n)\seq \A_1(n) \seq \A_2(n)$, 
$\A_0(1)=Z[\C[F_1]]$  and $\A_2(1)=\C[F_1]$.

For each $\lambda\in F_n^{\wedge}$ and $\mu\in F_1^{\wedge}$, let
$m_{\lambda\mu}$ be the number of GZ-subspaces of $V^{\lambda}$ isomorphic to
$V^{\mu}$, i.e., $m_{\lambda\mu}$ is the number of directed paths from $\mu$
to $\lambda$ in the branching graph. 
It is easily seen that
$\A_1(n)$ is a maximal commutative subalgebra of $\C[F_n]$ and that 
\beq \label{dim1}
\mbox{dim }\A_0(n) &=& 
\sum_{\lambda\in F_n^{\wedge}}\;\sum_{\mu\in F_1^{\wedge}}
m_{\lambda\mu},\\ \label{dim2}
\mbox{dim }\A_1(n) &=& \sum_{\lambda\in F_n^{\wedge}}\;\sum_{\mu\in F_1^{\wedge}}
m_{\lambda\mu}\;\mbox{dim }V^{\mu},\\ \label{dim3}
\mbox{dim }\A_2(n) &=& \sum_{\lambda\in F_n^{\wedge}}\;\sum_{\mu\in F_1^{\wedge}}
m_{\lambda\mu}\;(\mbox{dim }V^{\mu})^2.
\eeq

We denote $Z[\C[F_i]]$ by $Z_i$.
\bt We have 

\noi (i) $\A_0(n) = \langle Z_1, Z_2, \ldots ,Z_n \rangle$. 

\noi (ii) $\A_1(n) = \langle \A_1(1), Z_1, Z_2, \ldots ,Z_n \rangle$. 

\noi (iii) $\A_2(n) = \langle \C[F_1], Z_1, Z_2, \ldots ,Z_n \rangle$. 
\et
\pf (i) Consider the chain $T$ from (\ref{ch}) above. For $i=1,2,\ldots ,n$,
let
$p_{\lambda_{i}}\in Z_i$ denote the primitive central idempotent corresponding to the 
representation $\lambda_i \in F_i^{\wedge}$. Define 
$p_T\in \langle Z_1, Z_2, \ldots ,Z_n \rangle$ 
by
$$p_T = p_{\lambda_1}p_{\lambda_2}\cdots p_{\lambda_n}.$$
A little reflection shows that the image of $p_T$ under the isomorphism
(\ref{iso}) is $(f_{\mu}\;:\; \mu\in F_n^{\wedge})$, where $f_{\mu} =  
0$, if $\mu\not= \lambda$ and $f_{\lambda}$ is the projection on $V_T$ (with
respect to the decomposition  (\ref{decomp}) of $V^{\lambda}$). The result
follows since the primitive central idempotents corresponding to the irreducible
representations of a finite group 
form a basis of the center of the group algebra of the group.

\noi (ii) Note that $\C[F_1]$ commutes with  $Z_1, \ldots ,Z_n$. The result now
follows from part (i) and the isomorphism (\ref{iso}) with $n=1$. 

\noi (iii) Similar to part (ii). $\Box$

We call $\A_0(n)$ the {\em 
Gelfand-Tsetlin algebra (GZ-algebra)} of the 
multiplicity free chain of groups (\ref{c})
and denote it by $GZ_n$. 
Following {\bf\cite{p}} we call $\A_2(n)$ the {\em generalized
Gelfand-Tsetlin algebra}. 
By a {\em $GZ$-subspace of $F_n$} we mean a 
$GZ$-subspace of some irreducible representation $V^{\lambda}$ of $F_n,\;
\lambda\in F_n^{\wedge}$. 
By a {\em $GZ$-vector of $F_n$} we mean a vector in some 
$GZ$-subspace of some irreducible representation $V^{\lambda}$ of $F_n,\;
\lambda\in F_n^{\wedge}$. 
As an immediate consequence of
the theorem above we get the following result.

\bl \label{sc}
(i) Let $v\in V^{\lambda},\;\lambda\in F_n^{\wedge}$. If $v$ is an
eigenvector (for the action) of every element of $GZ_n$, then 
$v$ belongs to some GZ-subspace of $V^{\lambda}$.

\noi (ii) Let $v,u$ be two GZ-vectors of $F_n$. 
If $v$ and $u$ have the same eigenvalues
for every element of $GZ_n$, then $v$ and $u$ belong to the same GZ-subspace
of $V^{\lambda}$, for some $\lambda \in F_n^{\wedge}$.
\el
   
In Section 4 we define a multiplicity free chain of subgroups of $G_n$ and
consider the corresponding GZ-algebras.

\section{\large{Simplicity of branching and Young-Jucys-Murphy elements}}

Let $M$ be a complex finite dimensional semisimple algebra and let $N$ be a
semisimple subalgebra. 
Define the {\em relative commutant} of this pair to be the subalgebra
\beqn
Z(M,N)&=& \{m\in M \;|\; mn=nm \mbox{ for all }n\in N\},
\eeqn
consisting of all elements of $M$ that commute with $N$. 

The following result is well known. We include a proof for completeness.
\bt \label{mfc} 
Let $M$ be a complex finite dimensional semisimple algebra and let $N$ be a
semisimple subalgebra. 
Then $Z(M,N)$
is semisimple and the
following conditions are equivalent:

\noi
1. The restriction of any finite dimensional complex irreducible
representation of $M$ to $N$ is multiplicity free.

\noi
2. The relative commutant $Z(M,N)$ is commutative.
\et 
\pf By Wedderburn's theorem we may assume, without loss of generality, that
$M = M_1 \oplus \cdots \oplus M_k$, where each $M_i$ is a matrix algebra.  
We write elements of $M$ as $(m_1,\ldots ,m_k)$, where $m_i\in M_i$.
For $i=1,\ldots ,k$, let $N_i$ denote the image of $N$ under the natural
projection of $M$ onto $M_i$. Being the homomorphic image of a semisimple
algebra, $N_i$ itself is semisimple.

We have $Z(M,N)=Z(M_1,N_1)\oplus \cdots \oplus Z(M_k,N_k)$. By the double
centralizer theorem each $Z(M_i,N_i)$, and thus $Z(M,N)$, is semisimple. 
     
For $i=1,\ldots ,k$, let $V_i$ denote the $M$-submodule consisting of all
$(m_1,\ldots ,m_k)\in M$ with $m_j = 0$ for $ j\not= i$ and with all entries
of $m_i$ not in the first column equal to zero.
Note that
$V_1,\ldots ,V_k$ are all the distinct inequivalent irreducible
$M$-modules and that the  
decomposition of $V_i$ into irreducible $N$-modules is identical to
the decomposition of $V_i$ into irreducible $N_i$-modules.

It now follows from the double centralizer theorem that
$V_i$ is multiplicity free as a $N_i$-module, for all $i$ if and only if all
irreducible
modules of $Z(M_i,N_i)$ have dimension $1$, for all $i$ 
if and only if $Z(M_i,N_i)$ is  
abelian, for all $i$ if and only if $Z(M,N)$ is abelian. $\Box$  

Define the following subalgebras of $\C[H_{n,n}]$:

(i) For $2\leq m \leq n$, set $Z_{m,m-1,n}=Z[\C[H_{m,n}],\C[H_{m-1,n}]]$.

(ii) For $1\leq m \leq n$, set $Z_{m,n}=Z[\C[H_{m,n}]]$. 

In this section we show that the branching from $H_{m,n}$ to 
$H_{m-1,n}$, $2\leq m \leq n$ is multiplicity free and find a convenient set
of generators, the Young-Jucys-Murphy elements,  of the Gelfand-Tsetlin
algebra of the chain of groups, 
\beq \label{cog}
&H_{1,n}\seq H_{2,n}\seq \cdots \seq H_{n,n},&
\eeq
over the center $Z_{1,n}$ of the group algebra of the base group $H_{1,n}=G^n$.

For $i=1,\ldots ,n-1$ and $j=1,\ldots ,t$ define the following 
elements of $\C[H_{n,n}]$:
$$
Y_{i,j} = \sum_{\tau} \tau,$$
where the sum is over all elements
$\tau\in H_{n,n}$ satisfying the following properties: $\tau$ is a 
nontrivial $i$-cycle of type $j$ and $n$ does not belong
to the nontrivial cycle. 
Note that $Y_{1,1}=0$ and that all other $Y_{i,j}$ are
nonzero. We also set $Y_{n,j}=0$ for all $1\leq j\leq t$.

For $i=1,\ldots ,n$ and $j=1,\ldots ,t$ define the following 
elements of $\C[H_{n,n}]$:
$$Y'_{i,j} = \sum_{\tau} \tau,$$
where the sum is over all elements
$\tau\in H_{n,n}$ satisfying the following properties: $\tau$ is a 
nontrivial $i$-cycle of type $j$ and $n$ belongs
to the nontrivial cycle. 
Note that $Y'_{1,1}=0$ and that all other $Y'_{i,j}$ are
nonzero.

Define
\beqn
\pc'_n (G_*) &=& \{(\rho ,\lambda ,j) |\; \rho\in \pc_n(G_*), \lambda \in
\Po, j\in \{1,\ldots ,t\} \mbox{ with }\lambda\in \rho(C_j)\}.
\eeqn

For $(\rho,\lambda,j)\in \pc'_n(G_*)$ define $c_{(\rho,\lambda,j)}\in
\C[H_{n,n}]$ to be the sum
of all elements $\tau \in H_{n,n}$ satisfying 

\noi (i) $\mbox{type}(\tau) = \rho$.

\noi (ii) size of the cycle of $\tau$ containing $n$ is $\lambda$ and the
corresponding cycle product is of type $j$.

Note that, for $1\leq i \leq n,\;1\leq j\leq t$, $Y_{i,j}$ and $Y'_{i,j}$ 
are equal to $c_{(\rho,\lambda,j)}$, for suitable choice of $\rho$,
$\lambda$, and $j$.

\bl \label{yjm0} 
(i) $\{c_{(\rho,\lambda,j)} \;|\; (\rho,\lambda,j) \in \pc'_n(G_*)\}$ is a
basis of $Z_{n,n-1,n}$. It follows that
$$ \langle Y_{i,j}, Y'_{i,j}\;|\; 1\leq i \leq n,\;1\leq j \leq t \rangle \seq
Z_{n,n-1,n}.$$

\noi (ii) For $(\rho,\lambda, j)\in \pc'_n(G_*)$ we have
$$c_{(\rho,\lambda,j)}\in 
\langle Y_{i,p}, Y'_{i,p}\;|\; 1\leq p \leq t,\;1\leq i \leq k  \rangle,$$ 
where $k=\#\rho$. 

\noi (iii) $Z_{n,n-1,n} =
\langle Y_{i,j}, Y'_{i,j}\;|\; 1\leq j \leq t,\;1\leq i \leq n  \rangle.$ 

\noi (iv) $Z_{n-1,n} =
\langle Y_{i,j}, Y'_{1,j}\;|\; 1\leq j \leq t,\;1\leq i \leq n-1  \rangle.$ 
\el
\pf (i) Let $\tau, \tau' \in H_{n,n}$. Then, using the same argument that
characterized conjugacy in $G_n$ in Section 2, we can show
that $\tau = \sigma
\tau' \sigma^{-1}$ for some $\sigma \in H_{n-1,n}$ iff 
$\mbox{type}(\tau)=\mbox{ type}(\tau')$, the length of the cycle containing
$n$ is same in both $\tau$ and $\tau'$, and the cycle products of the cycles
containing $n$ are of the same type in both $\tau$ and $\tau'$. The result
follows.

\noi (ii) By induction on $\#\rho$. If $\#\rho = 0$, then $c_{(\rho,\lambda,j)}$
is the identity element of $\C[H_{n,n}]$ and the result is clearly true. Assume
the result whenever
$\#\rho \leq l$. Consider $(\rho,\lambda,j)\in \pc'_n(G_*)$ with
$\#\rho= l + 1$. 

Denote the multiset of nontrivial parts of $\rho$ by
$\{(\lambda_1,j_1),(\lambda_2,j_2),\ldots ,(\lambda_m,j_m)\}$.
Consider the following two subcases:

(a) $\lambda=1$ and $j=1$: Consider the product 
$Y_{\lambda_1,j_1}Y_{\lambda_2,j_2}\cdots Y_{\lambda_m,j_m}$.
Using part (i) we see that this product is in $Z_{n,n-1,n}$ and thus can be
expanded in the basis given in part (i). A little reflection shows that
\beqn
Y_{\lambda_1,j_1}Y_{\lambda_2,j_2}\cdots Y_{\lambda_m,j_m}&=&
\alpha_{(\rho,\lambda,j)}c_{(\rho,\lambda,j)} + 
\sum_{(\rho',\lambda',j')}
\alpha_{(\rho',\lambda',j')}c_{(\rho',\lambda',j')}, 
\eeqn
where $\alpha_{(\rho,\lambda,j)}\in \Po$, 
$\alpha_{(\rho',\lambda',j')}\in \N$ and 
the sum is over all $(\rho',\lambda',j')$ with $\#\rho' < \#\rho$. 
The result follows by induction.

(b) $\lambda\not=1$ or $j\not=1$: Without loss of generality we may assume
$(\lambda,j) = (\lambda_1,j_1)$. Now consider the product 
\beqn
Y'_{\lambda_1,j_1}Y_{\lambda_2,j_2}\cdots Y_{\lambda_m,j_m}&=&
\alpha_{(\rho,\lambda,j)}c_{(\rho,\lambda,j)} + 
\sum_{(\rho',\lambda',j')}
\alpha_{(\rho',\lambda',j')}c_{(\rho',\lambda',j')}, 
\eeqn
where $\alpha_{(\rho,\lambda,j)}\in \Po$, 
$\alpha_{(\rho',\lambda',j')}\in \N$ and 
the sum is over all $(\rho',\lambda',j')$ with $\#\rho' < \#\rho$. 
The result follows by induction.

\noi (iii) Follows from parts (i) and (ii).

\noi (iv) Embed $H_{n-1,n-1}$ into $H_{n-1,n}$ in the obvious way giving rise to
an embedding $\phi:Z_{n-1,n-1} \rar Z_{n-1,n}$. Note that $Z_{n-1,n}$ is
isomorphic to the tensor product of $\phi(Z_{n-1,n-1})$ and $Z[\C[G^{(n)}]]$.
Now $Z[\C[G^{(n)}]]$ is generated by $Y'_{1,j}, 1\leq j \leq t$ and a proof
similar to the proof of part (iii) shows that $\phi(Z_{n-1,n-1})$ is generated
by $Y_{i,j}, 1\leq j \leq t, 1\leq i \leq n-1$. The result follows. $\Box$

For $i=1,\ldots ,n$ define the following elements of $\C[H_{n,n}]$:
\beqn  X_i & = & \sum_{k=1}^{i-1} \sum_{g\in G} 
(g^{-1})^{(k)} g^{(i)} (k,i). \eeqn
Note that $X_1=0$. It is easy to see that 
$X_i$ is the difference of an element in $Z_{i,n}$
and an element in $Z_{i-1,n}$.
These elements are
called the {\em Young-Jucys-Murphy} (YJM) elements.

\bt \label{yjm1}
(i)
$Z_{m,m-1,n} = \langle Z_{m-1,n}, X_m \rangle ,\; 2\leq m \leq n.$

\noi (ii) $Z_{m,m-1,n},\; 2\leq m \leq n$ is commutative.

\et
\pf (i) We first consider $m=n$. Clearly $Z_{n,n-1,n} \supseteq \langle
Z_{n-1,n}, X_n \rangle$ (note that $X_n=Y'_{2,1}$).
To show the converse we need to show (by Lemma \ref{yjm0} (iii) and (iv)) 
that
$Y'_{i,j} \in \langle Z_{n-1,n}, X_n \rangle$, for $i=2,\ldots ,n$ and
$j=1,\ldots ,t$. 

Observe that, for $2\leq i\leq n$ and $2\leq j \leq t$, we have
\beq \label{ind1} Y'_{i,j}&=&Y'_{1,j}Y'_{i,1}.\eeq
Since $Y'_{1,j} \in Z_{n-1,n}$, for $1\leq j \leq t$ it is enough to show that
$Y'_{i,1} \in \langle Z_{n-1,n}, X_n \rangle$, for $i=2,\ldots ,n$. We show
this by induction on $i$.

Since $Y_{2,1}' = X_n$ we have
$Y_{2,1}' \in  \langle Z_{n-1,n}, X_n \rangle $.
Suppose $Y_{2,1}',\ldots , Y_{k+1,1}' \in  \langle Z_{n-1,n}, X_n \rangle $. We
shall now show that $Y_{k+2,1}'\in  \langle Z_{n-1,n}, X_n \rangle $.

We write $Y_{k+1,1}'$ as 
$$
\sum_{i_1,\ldots ,i_k,g_1,\ldots ,g_{k+1}} 
g_1^{(i_1)}g_2^{(i_2)}\cdots g_k^{(i_k)}g_{k+1}^{(n)}
(i_1,\ldots ,i_k , n),$$
where the sum is over all $(i_1,\ldots ,i_k) \in \{1,\ldots
,n-1\}^k$ with distinct components and all 
$(g_1,\ldots,g_{k+1})\in G^{k+1}$ with
$g_{k+1}g_k\cdots g_2g_1=e$.
In the following we use this summation convention implicitly.

Now consider the product $Y_{k+1,1}' X_n\in  \langle Z_{n-1,n}, X_n \rangle $:
\beq \label{yjm}&\ds{
\left\{\sum_{i_1,\ldots ,i_k,g_1,\ldots,g_{k+1}} 
g_1^{(i_1)}g_2^{(i_2)}\cdots g_k^{(i_k)}g_{k+1}^{(n)}
(i_1,\ldots ,i_k , n)\right\}
\left\{\sum_{i=1}^{n-1} \sum_g (g^{-1})^{(i)} g^{(n)} 
(i,n)\right\}.}& \eeq
Take a typical element 
$$g_1^{(i_1)}g_2^{(i_2)}\cdots g_k^{(i_k)}g_{k+1}^{(n)}
(i_1,\ldots ,i_k , n) (g^{-1})^{(i)} g^{(n)} (i,n)$$
of this product.
If $i\not= i_l$, for $l=1,\ldots ,k$, this product is 
$$(g^{-1})^{(i)}(g_1g)^{(i_1)}g_2^{(i_2)}\cdots g_k^{(i_k)}g_{k+1}^{(n)}
(i,i_1,\ldots ,i_k , n).$$
Note that $g_{k+1}\cdots g_2(g_1g)(g^{-1})=e$.

On the other hand if $i=i_l$, for some $1\leq l\leq k$, this product
becomes 
$$(g_1g)^{(i_1)}g_2^{(i_2)}\cdots
g_l^{(i_l)}(g_{l+1}g^{-1})^{(i_{l+1})}g_{l+2}^{(i_{l+2})}\cdots
g_{k+1}^{(n)}(i_1,\ldots ,i_l)(i_{l+1},\ldots ,n).$$
Note that, since $g_{k+1}g_k\cdots g_1=e$, we have 
$$g_{k+1}\cdots g_{l+2}(g_{l+1} g^{-1}) = g (g_l\cdots g_2 (g_1 g))^{-1} g^{-1}.$$

It follows that the element (\ref{yjm}) above is equal to
\beq \label{yjm2}
\lefteqn{\ds{\sum_{i,i_1,\ldots ,i_k,g,g_1,\ldots ,g_{k+1}} 
g^{(i)}g_1^{(i_1)}\cdots g_k^{(i_k)}g_{k+1}^{(n)}
(i,i_1,\ldots ,i_k ,n)}}\\ \nonumber
&&  + \ds{\sum_{i_1,\ldots ,i_k} \sum_{l=1}^k \sum_{j=1}^t
\sum_{g_1,\ldots ,g_{k+1}} \frac{|G|}{|C_j|}
g_1^{(i_1)}\cdots g_{k}^{(i_k)}g_{k+1}^{(n)}
(i_1,\ldots ,i_l)(i_{l+1},\ldots ,i_k , n),}
 \eeq
where the first sum is over all $(i,i_1,\ldots ,i_k) \in \{1,2,\ldots
,n-1\}^{k+1}$ with distinct components and all $(g,g_1,\ldots ,g_{k+1})\in G^{k+2}$ with
$g_{k+1}\cdots g_1 g = e$
and the second sum is over all $(i_1,\ldots ,i_k) \in
\{1,2,\ldots ,n-1\}^k$ with distinct components and all 
$(g_1,\ldots ,g_{k+1})\in G^{k+1}$
with type of $g_{k+1}\cdots g_{l+1}$ equal to $j$ and the type of 
$g_l\cdots g_1$ equal to $\I(j)$. 

We can rewrite (\ref{yjm2}) as
$$ Y_{k+2,1}' + \sum_{(\rho,\lambda,j)} \alpha_{(\rho,\lambda,j)} 
c_{(\rho,\lambda,j)},$$
where $\alpha_{(\rho,\lambda,j)}\in \N$ and the sum is over all $(\rho,\lambda,j)$ with $\#\rho \leq k+1$. By
induction hypothesis, (\ref{ind1}), 
and part (ii) of Lemma \ref{yjm0} it follows that
$Y_{k+2,1}'\in  \langle Z_{n-1,n}, X_n \rangle $.

We have now shown $Z_{n,n-1,n}=\langle Z_{n-1,n},X_n \rangle$. The case of
general $Z_{m,m-1,n}$ can be shown by embedding $Z_{m,m-1,m}$ in
$Z_{m,m-1,n}$ (as in part (iv) of Lemma \ref{yjm0}).

\noi (ii) This follows from part (i) since $X_m \in Z_{m,n} - Z_{m-1,n}$. $\Box$

It follows from Theorem \ref{mfc} and part (ii) of Theorem \ref{yjm1} that
the chain (\ref{cog}) is multiplicity free. Set, for $1\leq m
\leq n$,
$$GZ_{m,n} =\langle Z_{1,n}, Z_{2,n},\ldots , Z_{m,n} \rangle,$$
so that $GZ_{n,n}$ is the Gelfand-Tsetlin algebra of the chain (\ref{cog}).
Note that $X_i \in GZ_{i,n} \subseteq GZ_{n,n}$.

\bt \label {yjm3} We have
$$GZ_{m,n} = \langle Z[\C[G^n]], X_1,X_2,\ldots ,X_m \rangle ,\;1\leq m \leq
n. $$
\et
\pf The proof is by induction on $n$ and, for each $n$, by induction on
$m$. The cases $n=1,2$ are clear. Now consider general $n$.
The case $m=1$ is obvious. 
Assume we have proved that
$GZ_{m-1,n} = \langle Z[\C[G^n]], X_1,X_2,\ldots ,X_{m-1} \rangle $. 
It remains to show
that $GZ_{m,n} = \langle GZ_{m-1,n}, X_m \rangle$. The left hand side 
clearly contains the
right hand side so it suffices to check that the left hand side 
is contained in the right hand side. For this
it suffices to check that $Z_{m,n} \subseteq \langle GZ_{m-1,n}, X_m \rangle$.  
This follows from part (i) of Theorem \ref{yjm1}
since $Z_{m,n} \subseteq Z_{m,m-1,n}$. $\Box$

Let $V$ be a GZ-subspace of $G_n$. Then $V$ is an irreducible
$G^n$-module and is thus isomorphic to $\rho_1\otimes \cdots \otimes
\rho_n,\;\rho_i\in G^{\wedge}$ for all $i$. We call $\rho=(\rho_1,\ldots
,\rho_n)$ the {\em label} of $V$.  

Define
\beqn \alpha(V)&=&(\rho,\alpha_1,\ldots ,\alpha_n)\in \C^n ,\eeqn
where $\alpha_i = \mbox{eigenvalue of $X_i$ on $V$}$. We call $\alpha(v)$ the
{\em weight} of $V$ (note that $\alpha_1 = 0$ since $X_1 = 0$).  Define the
{\em spectrum} of $G_n$ by
$$ \spec(n) = \left\{ \alpha(V)\;:\; V \mbox{ is a GZ-subspace of }G_n
\right\}.$$
Let $V$ be a GZ-subspace of $G_n$ with label $\rho$. Then the primitive central
idempotent in $Z[\C[G^n]]$ corresponding to $\rho$ will have eigenvalue 1 on
$V$ and eigenvalue 0 on GZ-subspaces with different labels.
It now follows from Lemma \ref{sc} and Theorem \ref{yjm3} that a GZ-subspace
is uniquely determined by its weight.

By definition of GZ-subspaces and Lemma \ref{sc}, 
the set $\spec (n)$ is in natural bijection
with chains 
\beq \label{c1} 
T&=&\lambda_1\nearrow\lambda_2 \nearrow \cdots \nearrow \lambda_n,
\eeq
where $\lambda_i\in H_{i,n}^{\wedge},\;1\leq i \leq n$, in the Bratelli
diagram of (\ref{cog}).

Given $\alpha \in \spec (n)$ we denote by $V_\alpha$ the GZ-subspace with
weight $\alpha$ and by $T_\alpha$ the corresponding chain in the branching
graph. Similarly, given a chain $T$ as in (\ref{c1}) we denote the
correponding GZ-subspace by $V_T$ and the 
weight vector $\alpha(V_T)$ by $\alpha(T)$. Thus we have 1-1
correspondences
$$ T \mapsto \alpha(V_T),\;\;\; \alpha \mapsto T_\alpha$$
between chains (\ref{c1}) and $\spec (n)$. For $\lambda\in
H_{n,n}^{\wedge}$ define
$$\spec(n,\lambda)=\{\alpha \in \spec(n) | T_{\alpha} \mbox{ ends at
}\lambda\}.
$$

We have, from (\ref{dim1}),
\beqn \mbox{dim }GZ_{n,n} &=& |\spec(n)|,\\
  \mbox{dim }V^{\lambda} &=& \sum_{\alpha\in\spec(n,\lambda)} \mbox{dim } V_{\alpha},\;\lambda
\in H_{n,n}^{\wedge}.
\eeqn

There is a natural equivalence relation $\sim$ on $\spec (n)$: for $\alpha,
\beta \in \spec (n)$,
\beqn \alpha \sim \beta &\Leftrightarrow& \alpha, \beta \in \spec(n,\lambda)
\mbox{ for some }\lambda \in H_{n,n}^{\wedge}.
\eeqn
Clearly we have $|\spec (n) / \sim| = |H_{n,n}^{\wedge}|$.

\section{\large{Action of Coxeter generators on GZ-subspaces}} 

In this section
we describe the action of the Coxeter generators on
GZ-subspaces in terms of transformations of weights. 

Let $\lambda\in H_{n,n}^{\wedge}$. We have the GZ-decomposition
\beq \label{gzdi}
V^{\lambda}&=&\oplus_{\alpha\in \spec(n,\lambda)} V_{\alpha},
\eeq
of $V^{\lambda}$ into irreducible $G^n$-modules.

We now consider the action of the Coxeter generators $s_i=(i,i+1)$ of $S_n$
on $V^{\lambda}$. Since the $V_{\alpha}$ consist of common eigenvectors of
$X_1,\ldots ,X_n$ and are $G^n$-invariant, it is useful to know the
commutation relations satisfied by the $s_i$, the $X_j$, and the $g^{(l)}$.

\bl \label{cr}
The following relations hold in $G_n$:

(i) $X_1,\ldots ,X_n$ commute.

(ii) $X_ig^{(l)}=g^{(l)}X_i$, $g\in G,\; 1\leq i,l \leq n$.

(iii) $s_ig^{(i)}s_i = g^{(i+1)}$, $g\in G,\; 1\leq i \leq n-1$. In particular,
$s_i^2=1$, $1\leq i\leq n-1$.

(iv) $s_ig^{(l)}=g^{(l)}s_i$, $1\leq i\leq n-1,\; 1\leq l \leq n,\; l\not= i,i+1$.

(v) $s_iX_is_i + \sum_{g\in G} g^{(i+1)}s_i(g^{-1})^{(i+1)}=X_{i+1}$,
$1\leq i \leq n-1$.

(vi)  $s_iX_l=X_ls_i$,  $1\leq i\leq n-1,\; 1\leq l \leq n,\; l\not= i,i+1$.

\el
\pf (i) We have already seen this.

\noi (ii) This can be checked directly. An alternative proof is as follows. On
every GZ-subspace of $G_n$, the actions of $X_i$ and $g^{(l)}$ clearly commute.
By considering the isomorphism 
\beqn
\C[G_n]&\cong&\bigoplus_{\lambda\in G_n^{\wedge}} \mbox{End}(V^{\lambda}),
\eeqn
given by
$$g \mapsto ( V^{\lambda} \buildrel {g}\over \rightarrow V^{\lambda}\;:\;
\lambda
\in G_n^{\wedge}),\;\;g\in G_n,$$
we see that $X_i$ and $g^{(l)}$ commute in $G_n$.

\noi (iii) and (iv) This is clear.

\noi (v) We have 
\beqn
s_iX_is_i &=& 
(i, i+1) \left( \sum_{k=1}^{i-1} \sum_{g\in G} (g^{-1})^{(k)} g^{(i)} (k, i)
\right) (i, i+1) \\
&=& \sum_{k=1}^{i-1} \sum_{g\in G} (g^{-1})^{(k)} g^{(i+1)} (k, i+1) \\
&=& X_{i+1} - \sum_{g\in G} (g^{-1})^{(i)} g^{(i+1)} (i, i+1).
\eeqn 

\noi (vi) First assume $l\leq i-1$. Then
\beqn
s_i X_l &=& 
(i, i+1) \left( \sum_{k=1}^{l-1} \sum_{g\in G}
(g^{-1})^{(k)}g^{(l)} (k,l)\right)\\ 
        &=&\left( \sum_{k=1}^{l-1} \sum_{g\in G}
(g^{-1})^{(k)}g^{(l)} (k,l)\right)(i, i+1)\\ 
        &=& X_l s_i.
\eeqn
Now assume $l\geq i+2$. Then
\beqn
s_i X_l 
&=& (i, i+1) \left( 
\sum_{k=1}^{i-1} \sum_{g\in G} (g^{-1})^{(k)} (g)^{(l)} (k, l)
+ \sum_{k=i+2}^{l-1} \sum_{g\in G} (g^{-1})^{(k)} (g)^{(l)} (k, l)\right. \\
&& + \left. \sum_{g\in G} (g^{-1})^{(i)} (g)^{(l)} (i, l)
+ \sum_{g\in G} (g^{-1})^{(i+1)} (g)^{(l)} (i+1, l)\right)\\
&=& \left( 
\sum_{k=1}^{i-1} \sum_{g\in G} (g^{-1})^{(k)} (g)^{(l)} (k, l)
+ \sum_{k=i+2}^{l-1} \sum_{g\in G} (g^{-1})^{(k)} (g)^{(l)} (k, l)\right. \\
&& + \left. \sum_{g\in G} (g^{-1})^{(i+1)} (g)^{(l)} (i+1, l)
+ \sum_{g\in G} (g^{-1})^{(i)} (g)^{(l)} (i, l)\right) (i, i+1)\\
&=& X_l s_i. \;\;\Box
\eeqn

Using part (iii) of Lemma \ref{cr} we can
rewrite part (v) of Lemma \ref{cr} as
\beq \label{ccr}
X_is_i + \sum_{g\in G} g^{(i)} (g^{-1})^{(i+1)} &=& s_iX_{i+1},\; 1\leq i \leq
n-1. 
\eeq

Consider the irreducible $G^n$-module $V_{\alpha}$ in the decomposition
(\ref{gzdi}) above. Let $V$ be the subspace of $V^{\lambda}$ spanned by
$V_{\alpha}$ and $s_i\cdot V_{\alpha}$. Lemma \ref{cr} shows that $V$ is
invariant under the actions of $s_i, X_i, X_{i+1}$, and $G^n$. A
study of this action will enable us to write down matrices for the action of
$s_i$ on $V_{\alpha}$. 

\bl \label{ss}
For $i=1,2,\ldots ,n-1$, let $A_i$ be the subalgebra of  
$\C[G_n]$ generated by 
$G^n, s_i, X_i$, and $X_{i+1}$. 
Then $A_i$ is semisimple and the actions of $X_i$
and $X_{i+1}$ on any $A_i$-module are simultaneously diagonalizable.
\el
\pf Let $Mat(n)$ denote the algebra of complex $(|G|^n n!) \times (|G|^n
n!)$ matrices, with rows and columns indexed by elements of $G_n$. 
Consider
the left regular representation of $G_n$. 
Writing this in matrix terms gives
an embedding of $\C[G_n]$ into $Mat(n)$. 
We write $\gamma :\C[G_n] \hookrightarrow Mat (n)$.

Note that 

(a) The left action of $(i,i+1)$ on $G_n$ is inverse to itself.

(b) For $k<i$, the left action of $(g^{-1})^{(k)} g^{(i)}(k,i)$
on $G_n$ is inverse to itself.

(c) For $g\in G$ and $1\leq l \leq n$, 
the left action of $g^{(l)}$ on $G_n$
is inverse to the action of $(g^{-1})^{(l)}.$

It follows that the matrices $\gamma(s_i), \gamma(X_i), \gamma(X_{i+1})$
are real and symmetric and that the generating set
$$\{\gamma(s_i), \gamma(X_i), \gamma(X_{i+1}) \} \cup \{\gamma(g^{(l)}) : 
1\leq l\leq n, g\in G\}$$
of $\gamma(A_i)$ is closed under the matrix $*$ operation
$M\mapsto (\bar M)^t$. So $\gamma(A_i)$ itself is closed 
under the $*$ operation and a standard    
result on finite dimensional $C^*$-algebras now shows that $\gamma(A_i)$
(and hence $A_i$) is semisimple. 

Part (ii) follows since  $\gamma(X_i)$ and $\gamma(X_{i+1})$ are
commuting real, symmetric matrices and thus the $*$-subalgebra of
$\gamma(A_i)$
generated by them is commutative. $\Box$

Before proceeding further we introduce some useful notation. 
For $1\leq i \leq n-1$, let $\omega_i$ be the involution on $\{1,2,\ldots
,n\}$ defined by $\omega_i(l)=l$, if $l\not=i,i+1$, $\omega_i(i)=i+1$, and
$\omega_i(i+1)=i$. Parts (iii) and (iv) of Lemma \ref{cr} may be written as
follows. For $1\leq i \leq n-1$ and $1\leq l \leq n$ we have
\beq \label{t1}
g^{(l)}s_i&=&s_ig^{\omega_i(l)}.
\eeq

Let $W_1$ and $W_2$ be vector spaces and set $W=W_1\otimes W_2$. We can
define a switch operator on $W$ that sends $w_1\otimes w_2$ to $w_2\otimes
w_1$. Now let $U$ be a vector space having the same dimension as $W$. We can
fix an isomorphism between $U$ and $W$ and transfer the switch operator on
$W$ to $U$ via this isomorphism. However, since the isomorphism between $U$
and $W$ is not canonical, there is no canonically defined switch operator on
$U$. This situation does not arise when we consider irreducible
$G^n$-modules.

Let $\rho=(\rho_1,\ldots ,\rho_n)$, where $\rho_i \in G^{\wedge}$ for all
$i$ and consider the irreducible $G^n$-module 
$V^{\rho_1}\otimes \cdots \otimes V^{\rho_n}$. 
Let $1\leq i \leq n-1$.
Define the involution, called the {\em switch operator}, 
$$
\tau_{i,\rho}: V^{\rho_1}\otimes \cdots \otimes V^{\rho_n} \rar 
V^{\rho_1}\otimes \cdots \otimes V^{\rho_{i-1}}
\otimes V^{\rho_{i+1}}\otimes V^{\rho_i}\otimes V^{\rho_{i+2}}\otimes \cdots
\otimes V^{\rho_n}$$
by switching the $i$ and $i+1$ factors: 
\beqn 
\tau_{i,\rho}(v_1\otimes \cdots \otimes v_n)&=&v_1\otimes \cdots \otimes
v_{i-1}\otimes v_{i+1}\otimes v_i \otimes v_{i+2} \otimes \cdots \otimes
v_n.
\eeqn
We have, for $g\in G,\;v\in
V^{\rho_1}\otimes \cdots \otimes V^{\rho_n},\; 1\leq l \leq n$,
\beq \label{t2}
\tau_{i,\rho}(g^{(l)}v)&=& g^{(\omega_i(l))}\tau_{i,\rho}(v).
\eeq
Now 
let $V$ be an irreducible $G^n$-module isomorphic to 
$V^{\rho_1}\otimes \cdots \otimes V^{\rho_n}$. 
Fix a $G^n$-linear  isomorphism 
$f: V\rar V^{\rho_1}\otimes \cdots \otimes V^{\rho_n}$. 
Define an involution 
$$\tau_{i,V}:V\rar V$$
by $\tau_{i,V}=f^{-1}\tau_{i,\rho} f$. 
It is easily seen by Schur's lemma that
$\tau_{i,V}$ is independent of the chosen $f$ and therefore $\tau_{i,V}$ is
canonically defined.
We have, for $g\in G,\;v\in V,\;1\leq l\leq n$,
\beq \label{t3} 
\tau_{i,V}(g^{(l)}v) &=& g^{(\omega_i(l))}
\tau_{i,V}(v).
\eeq
In what follows we shall use (\ref{t1}), (\ref{t2}), (\ref{t3}) (and
(\ref{t4}), (\ref{t5}), (\ref{t6}), (\ref{bvne}), (\ref{bve}) below) 
without explicit mention.

Let $V_{\alpha}$ be a GZ-subspace of $V^{\lambda}, \lambda\in
H_{n,n}^{\wedge}$ with weight 
$\alpha=(\rho,\alpha_1,\ldots ,\alpha_n)$, where $\rho=(\rho_1,\ldots
,\rho_n)$. 
Fix a $G^n$-linear isomorphism
$$f: V_{\alpha} \rar V^{\rho_1}\otimes \cdots \otimes V^{\rho_n}.$$
Let $1\leq i \leq n-1$.
Since $s_i^2=1$ the map $v\mapsto s_i\cdot v$ on $V^{\lambda}$ 
is an involution. 
Consider the subspace $s_iV_{\alpha}$ of $V^{\lambda}$. 
Then, by Lemma \ref{cr} (iii) and (iv),
$s_iV_{\alpha}$ is also closed under the $G^n$-action.
The map
\beq \label{t}  
&f^{\tau_i}: s_iV_{\alpha} \rar 
V^{\rho_1}\otimes \cdots \otimes V^{\rho_{i-1}}
\otimes V^{\rho_{i+1}}\otimes V^{\rho_i}\otimes V^{\rho_{i+2}}\otimes \cdots
\otimes V^{\rho_n},
\eeq
given by $f^{\tau_i}(s_iv)=\tau_{i,\rho}(f(v))$ 
is a $G^n$-linear isomorphism. 
To see this, let
$v\in V_{\alpha}$. 
Then, for $1\leq l \leq n$, we have
\beqn
 g^{(l)}f^{\tau_i}(s_iv)&=&
 g^{(l)}(\tau_{i,\rho}(f(v))),\\
f^{\tau_i}(g^{(l)}s_iv)&=&
f^{\tau_i}(s_i g^{(\omega_i(l))}v))\\
&=&\tau_{i,\rho}(f(g^{(\omega_i(l))}v))\\
&=&\tau_{i,\rho}(g^{(\omega_i(l))}(f(v)))\\
&=&g^{(l)}(\tau_{i,\rho}(f(v))).
\eeqn

For $1\leq i\leq n-1$, 
define an element $b_i=\sum_{g\in G} g^{(i)}(g^{-1})^{(i+1)} \in \C[G_n]$. 
For $h\in G$ we have $h^{(l)}b_i=b_ih^{(l)},\;l\not=i,i+1$ and
\beqn
&h^{(i)}b_i=\sum_{g\in G}(hg)^{(i)}(g^{-1})^{(i+1)}
=\sum_{g\in G}(hg)^{(i)}((hg)^{-1})^{(i+1)}h^{(i+1)}=b_ih^{(i+1)}.&
\eeqn
Similarly, we can show $h^{(i+1)}b_i=b_ih^{(i)}$. We have, for $1\leq l\leq
n$ and $h\in G$,
\beq \label{t4}
h^{(l)}b_i &=& b_i h^{(\omega_i(l))}. 
\eeq
Also note that we can rewrite (\ref{ccr}) as follows
\beq \label{t5}
X_is_i&=& s_iX_{i+1} - b_i,\;\;X_{i+1}s_i = s_iX_i + b_i.
\eeq

The map $V_{\alpha}\rar s_iV_{\alpha}$ given by 
$v\mapsto s_ib_iv$ is a $G^n$-linear map. 
This follows from, for $1\leq l \leq n$,
\beq \label{t6}
s_ib_i(g^{(l)}v)&=&s_i(g^{(\omega_i(l))}b_iv)\; = \;g^{(l)}s_ib_iv.
\eeq

It follows that 
\beq \label{bvne}
\rho_i\not= \rho_{i+1} &\mbox{implies}& b_iv=0,\;v\in V_{\alpha}.
\eeq
Now assume $\rho_i=\rho_{i+1}$. 
The map $V_{\alpha}\rar s_iV_{\alpha}$ given by $v\mapsto
s_i\tau_{i,V_{\alpha}}(v)$ is a $G^n$-linear
isomorphism. 
This follows from, for $1\leq l \leq n$,
\beqn
s_i\tau_{i,V_{\alpha}}(g^{(l)}v)&=&
s_i(g^{(\omega_i(l))}\tau_{i,V_{\alpha}}(v)) \;=\;
g^{(l)}s_i\tau_{i,V_{\alpha}}(v).
\eeqn
It follows that 
$b_iv=\beta\tau_{i,V_{\alpha}}(v), v\in
V_{\alpha}$,
for some scalar $\beta$. Now the trace of the action of $b_i$ on
$V^{\rho_1}\otimes \cdots \otimes V^{\rho_n}$ is 
$$
\frac{\dim(V^{\rho_1})\cdots
\dim(V^{\rho_n})}{\dim(V^{\rho_i})\dim(V^{\rho_{i+1}})}
\sum_{g\in G}\chi(g)\chi(g^{-1})=
\frac{\dim(V^{\rho_1})\cdots
\dim(V^{\rho_n})}{\dim(V^{\rho_i})\dim(V^{\rho_{i+1}})}|G|,
$$
by
the first orthogonality relation for characters ($\chi = $ character of
$V^{\rho_i})$. 
Since
the trace of $\tau_{i,V_{\alpha}}$ is 
$$\frac{\dim(V^{\rho_1})\cdots
\dim(V^{\rho_n})}{\dim(V^{\rho_i})},$$
it follows that $\beta= \frac{|G|}{\mbox{dim}(V^{\rho_i})}$. 
We have 
\beq \label{bve}
\rho_i = \rho_{i+1} &\mbox{implies}&
b_iv=\frac{|G|}{\mbox{dim}(V^{\rho_i})}\tau_{i,V_{\alpha}}(v),\;v\in
V_{\alpha}.
\eeq

The following result relates the action of $s_i$ on GZ-subspaces to
transformations on the corresponding weights.
\bt\label{mtl}
Let $\alpha=((\rho_1,\ldots ,\rho_n),\alpha_1,\ldots ,\alpha_n)\in
\spec(n,\lambda)$   
and consider the GZ-subspace $V_{\alpha}$ of $V^{\lambda}$. Then

\noi (i) For $1\leq i \leq n-1$, $s_i\cdot V_{\alpha} = V_{\alpha}$ iff 
$\rho_i=\rho_{i+1}$ and $\alpha_{i+1}=\alpha_i \pm
\frac{|G|}{\mbox{dim}(V^{\rho_i})}$.

\noi (ii) For $1\leq i \leq n-1$ we have 

(a) $\rho_i=\rho_{i+1}$ and $\alpha_{i+1}=\alpha_i +
\frac{|G|}{\mbox{dim}(V^{\rho_i})}$ implies $s_iv=\tau_{i,V_{\alpha}}(v),\;v\in
V_{\alpha}$.

(b) $\rho_i=\rho_{i+1}$ and $\alpha_{i+1}=\alpha_i -
\frac{|G|}{\mbox{dim}(V^{\rho_i})}$ implies $s_iv=-\tau_{i,V_{\alpha}}(v),\;v\in
V_{\alpha}$.

\noi (iii) If $\rho_i=\rho_{i+1}$ then $\alpha_i \not= \alpha_{i+1}$, $1\leq i \leq
n-1$.

\noi (iv) For $i=1,\ldots ,n-2$ the following statements are not true.

 (a) $\rho_i=\rho_{i+1}=\rho_{i+2}$ and $\alpha_i=\alpha_{i+1} +
\frac{|G|}{\mbox{dim}(V^{\rho_i})} = \alpha_{i+2}$.

 (b) $\rho_i=\rho_{i+1}=\rho_{i+2}$ and $\alpha_i=\alpha_{i+1} -
\frac{|G|}{\mbox{dim}(V^{\rho_i})} = \alpha_{i+2}$.

\noi (v) For $1\leq i \leq n-1$, if $\rho_i\not= \rho_{i+1}$ then $U=s_i\cdot
V_{\alpha}$ is a GZ-subspace of $V^{\lambda}$ with weight 
$$s_i\cdot \alpha =
((\rho_1,\ldots ,\rho_{i-1},\rho_{i+1},\rho_i,\rho_{i+2},\ldots ,\rho_n),
\alpha_1,\ldots ,\alpha_{i-1},\alpha_{i+1},\alpha_i,\alpha_{i+2},\ldots
,\alpha_n).
$$

\noi (vi) For $1\leq i \leq n-1$, if $\rho_i=\rho_{i+1}$ and
$\alpha_{i+1}\not= \alpha_i \pm \frac{|G|}{\mbox{dim}(V^{\rho_i})}$ then,
setting
$$U=\left(s_i -
\frac{|G|}
{(\alpha_{i+1}-\alpha_i)\mbox{dim}
(V^{\rho_i})}\tau_{i,V_{\alpha}}\right)V_{\alpha},
$$
we have that $U$ is a GZ-subspace $V^{\lambda}$ with weight  
$$s_i\cdot \alpha =
((\rho_1,\ldots ,\rho_n),
\alpha_1,\ldots ,\alpha_{i-1},\alpha_{i+1},\alpha_i,\alpha_{i+2},\ldots
,\alpha_n).
$$
\et
\pf (i) (only if) That $\rho_i=\rho_{i+1}$ is clear from (\ref{t}). 
The maps $V_{\alpha}\rar V_{\alpha}$ given by $v\mapsto
s_iv$ and $v\mapsto \tau_{i,V_{\alpha}}(v)$ are both involutions. 
Therefore the possible
eigenvalues are 1 and -1.

Let $u,v\in V_{\alpha}, u,v\not=0$ with $s_iv=v$ and $s_iu=-u$. Then
\beq \label{f}
&\alpha_i v=X_i(v)=X_i(s_iv)=(s_iX_{i+1}-b_i)v=
\alpha_{i+1} v-\frac{|G|}{\mbox{dim}(V^{\rho_i})}\tau_{i,V_{\alpha}}(v),&\\ 
\label{s}
&\alpha_i u=X_i(u)=X_i(-s_iu)=-(s_iX_{i+1}-b_i)u=
\alpha_{i+1} u+\frac{|G|}{\mbox{dim}(V^{\rho_i})}\tau_{i,V_{\alpha}}(u).& 
\eeq
It follows that $\tau_{i,V_{\alpha}}(v)$ 
is a multiple of $v$ and is thus either $v$ or
$-v$.
Similarly, $\tau_{i,V_{\alpha}}(u)$ is a multiple of $u$ 
and is thus either $u$ or
$-u$. Since there exists an eigenvector of $s_i:V_{\alpha}\rar V_{\alpha}$, 
the result follows.

(if) Since $s_iV_{\alpha}$ is also an irreducible $G^n$-module 
(by (\ref{t})) either
$V_{\alpha}=s_iV_{\alpha}$ or $s_iV_{\alpha}\cap V_{\alpha} =\{0\}$. 
Assume that $s_iV_{\alpha}\cap V_{\alpha} = \{0\}$. We shall derive
a contradiction.

We assume that 
$\alpha_{i+1}=\alpha_i + \frac{|G|}{\mbox{dim}(V^{\rho_i})}$. The case
$\alpha_{i+1}=\alpha_i - \frac{|G|}{\mbox{dim}(V^{\rho_i})}$ is similar.

The subspace $V_{\alpha}\oplus s_iV_{\alpha}$ of $V^{\lambda}$ 
is an $A_i$-module, by Lemma \ref{cr}.
Define a subspace
\beqn
M&=& \{ v - s_i\cdot \tau_{i,V_{\alpha}}(v) | v\in V_{\alpha}\} 
\seq V_{\alpha} \oplus s_iV_{\alpha}.
\eeqn
We check that $M$ is an $A_i$-submodule, i.e., is closed under the action of
$s_i,X_i,X_{i+1}$, and $g^{(l)}, l=1,\ldots ,n$. 
We have, for $v\in V_{\alpha}$,
\beqn
s_i(v-s_i\tau_{i,V_{\alpha}}(v))&=&s_iv-\tau_{i,V_{\alpha}}(v)\\
   &=&-\tau_{i,V_{\alpha}}(v)
-s_i(\tau_{i,V_{\alpha}}(-\tau_{i,V_{\alpha}}(v)))\in M, \\
g^{(l)}(v-s_i\tau_{i,V_{\alpha}}(v))&=&
g^{(l)}v-g^{(l)}s_i\tau_{i,V_{\alpha}}(v)\\
&=&
g^{(l)}v-s_ig^{(\omega_i(l))}\tau_{i,V_{\alpha}}(v)\\
&=& g^{(l)}v-s_i\tau_{i,V_{\alpha}}(g^{(l)}v)\in M, \\
X_i(v-s_i\tau_{i,V_{\alpha}}(v))&=&\alpha_iv-(s_iX_{i+1}-b_i)
(\tau_{i,V_{\alpha}}(v))\\
&=& \alpha_iv-\alpha_{i+1}s_i\tau_{i,V_{\alpha}}(v) 
+\frac{|G|}{\mbox{dim}(V^{\rho_i})}v\\
&=&
\alpha_{i+1}(v-s_i\tau_{i,V_{\alpha}}(v))\in M, \\
X_{i+1}(v-s_i\tau_{i,V_{\alpha}}(v))&=&\alpha_{i+1}v-(s_iX_i+b_i)
(\tau_{i,V_{\alpha}}(v))\\
&=& \alpha_{i+1}v-\alpha_is_i\tau_{i,V_{\alpha}}(v) 
-\frac{|G|}{\mbox{dim}(V^{\rho_i})}v\\
&=&
\alpha_i(v-s_i\tau_{i,V_{\alpha}}(v))\in M.
\eeqn

We shall now show that $M$ is the only nonempty, proper $A_i$-submodule of
$V_{\alpha}\oplus s_iV_{\alpha}$. Since 
$\mbox{dim}(M)<\mbox{dim}(V_{\alpha}\oplus s_iV_{\alpha})$, this contradicts
the fact that $V_{\alpha}\oplus s_iV_{\alpha}$ is a semisimple $A_i$-module
(since $A_i$ is a semisimple algebra by Lemma \ref{ss}). 

Let $M'$ be a nonempty, proper $A_i$-submodule of $V_{\alpha}\oplus
s_iV_{\alpha}$. Since $M'$ is
closed under $s_i$ we have $M'\not\seq V_{\alpha}$ and $M'\not\seq
s_iV_{\alpha}$. 
Also, $M'$ is
in particular a $G^n$-submodule of $V_{\alpha}\oplus s_iV_{\alpha}$. 
Since $V_{\alpha}$ and $s_iV_{\alpha}$
are isomorphic irreducible $G^n$-modules and 
$v\mapsto s_i\tau_{i,V_{\alpha}}(v)$ is
a $G^n$-linear isomorphism between them, it follows by Schur's lemma
that
\beqn
M'&=& \{v+\gamma s_i\tau_{i,V_{\alpha}}(v) | v\in V_{\alpha}\},
\eeqn
for some $0\not=\gamma \in \C$. We shall show that $\gamma = -1$. 

Now
$$s_i(v+\gamma s_i\tau_{i,V_{\alpha}}(v))=s_iv+\gamma\tau_{i,V_{\alpha}}(v)
=\gamma\tau_{i,V_{\alpha}}(v)
+\frac{\gamma s_i\tau_{i,V_{\alpha}}
(\gamma\tau_{i,V_{\alpha}}(v))}{\gamma^2},$$
which yields $\gamma^2=1$ and hence $\gamma = \pm 1$.

We have
\beqn
X_i(v+s_i\tau_{i,V_{\alpha}}(v))&=&\alpha_iv +
\alpha_{i+1}s_i\tau_{i,V_{\alpha}}(v) - b_i\tau_{i,V_{\alpha}}(v)\\
&=&
\left(\alpha_i-\frac{|G|}{\mbox{dim}(V^{\rho_i})}\right)v +
\alpha_{i+1}s_i\tau_{i,V_{\alpha}}(v),\\
X_i(v-s_i\tau_{i,V_{\alpha}}(v))&=&\alpha_iv -
\alpha_{i+1}s_i\tau_{i,V_{\alpha}}(v) + b_i\tau_{i,V_{\alpha}}(v)\\
&=&
\left(\alpha_i+\frac{|G|}{\mbox{dim}(V^{\rho_i})}\right)v -
\alpha_{i+1}s_i\tau_{i,V_{\alpha}}(v).
\eeqn
and thus $\gamma = -1$. Thus $M'=M$ and the proof of the if part is
complete.

\noi (ii) This follows from (\ref{f}) and (\ref{s}).

\noi (iii) Either $s_iV_{\alpha}=V_{\alpha}$ or $s_iV_{\alpha}\cap V_{\alpha} =\{0\}$.
If $s_iV_{\alpha}=V_{\alpha}$ then by part (i) $\alpha_{i+1}=\alpha_i \pm
\frac{|G|}{\mbox{dim}(V^{\rho_i})}$, so $\alpha_i\not= \alpha_{i+1}$.

Now assume $s_iV_{\alpha}\cap V_{\alpha}=\{0\}$. 
Then, as checked before, $V_{\alpha}\oplus s_iV_{\alpha}$ is
$A_i$-invariant. Choose a basis $B$ of $V_{\alpha}$ 
and consider the 
basis $B\cup s_iB$ of $V_{\alpha}\oplus s_iV_{\alpha}$. 
Let $N$ be the matrix of $\tau_{i,V_{\alpha}}$ with
respect to the basis $B$ and set $\kappa =\frac{|G|}{ \mbox{dim}(V^{\rho_i})}$.
Using the relation $X_is_i=s_iX_{i+1}-b_i$ we see that the
matrices of $X_i$ and $X_{i+1}$ (respectively) 
with respect to $B\cup sB$ are given in block
form as follows
$$
\left[ \ba{rr} \alpha_i I & -\kappa N\\
                  0       & \alpha_{i+1} I
        \ea
\right],\;\;
\left[ \ba{rr} \alpha_{i+1} I & \kappa N\\
                  0       & \alpha_i I
        \ea
\right]
$$ 
The actions of $X_i$ and $X_{i+1}$ on $V_{\alpha}\oplus s_iV_{\alpha}$ 
are diagonalizable by Lemma \ref{ss}
and thus $\alpha_i\not=\alpha_{i+1}$ (since $N\not=0$).

\noi (iv) Suppose statement (a) is true. Since, as a $G^n$-module, $V_{\alpha}$
is isomorphic to $V^{\rho_1}\otimes \cdots \otimes V^{\rho_n}$, we can
choose a $v\in V_{\alpha}$
such that $\tau_{i,V_{\alpha}}(v)=\tau_{i+1,V_{\alpha}}(v)=v$. By part (ii)
we have $s_iv=-v$ and $s_{i+1}v=v$. Now consider the Coxeter relation
$$s_is_{i+1}s_i=s_{i+1}s_is_{i+1}$$
and let both sides act on $v$. The left hand side yields $v$ and the right
hand side yields $-v$, a contradiction. So, statement (a) must be false. The
proof for the falsity of statement (b) is similar.

\noi (v) When $\rho_i\not=\rho_{i+1}$ then, by (\ref{bvne}), 
$b_iv=0$ for $v\in V_{\alpha}$. The result now
follows from (\ref{t}) and the relation $X_is_i=s_iX_{i+1}-b_i$. 

\noi (vi) By part (i) $s_iV_{\alpha}\cap V_{\alpha}=\{0\}$ and 
by part (iii) $\alpha_i\not=\alpha_{i+1}$. 
Clearly, $U$ is a subspace of $V_{\alpha}\oplus s_iV_{\alpha}$ and
since $v\mapsto s_i\tau_{i,V_{\alpha}}(v)$ 
(or, equivalently, $\tau_{i,V_{\alpha}}(v)\mapsto s_iv$)
is a $G^n$-linear isomorphism between $V_{\alpha}$
and $s_iV_{\alpha}$ it follows that $U$ is also an 
irreducible $G^n$-module with
label $(\rho_1,\ldots ,\rho_n)$. 
It remains to check that $X_i, X_{i+1}$ act on $U$  
by appropriate scalars. Setting $\kappa = \frac{|G|}{\mbox{dim}(V^{\rho_i})}$,
we have, using (\ref{bve}),
\beqn
X_i\left(s_iv-\frac{\kappa}{\alpha_{i+1}-\alpha_i}\tau_{i,V_{\alpha}}(v)\right)
&=& \alpha_{i+1}(s_iv)-b_iv
-\frac{\alpha_i\kappa}{\alpha_{i+1}-\alpha_i}\tau_{i,V_{\alpha}}(v)\\
&=&\alpha_{i+1}\left(s_iv-\frac{\kappa}{\alpha_{i+1}-\alpha_i}
\tau_{i,V_{\alpha}}(v)\right),\\
X_{i+1}\left(s_iv-\frac{\kappa}{\alpha_{i+1}-\alpha_i}\tau_{i,V_{\alpha}}(v)\right)
&=& \alpha_i(s_iv)+b_iv
-\frac{\alpha_{i+1}\kappa}{\alpha_{i+1}-\alpha_i}\tau_{i,V_{\alpha}}(v)\\
&=&\alpha_i\left(s_iv-\frac{\kappa}{\alpha_{i+1}-\alpha_i}
\tau_{i,V_{\alpha}}(v)\right).
\eeqn

That completes the proof. $\Box$

Let $\alpha = ((\rho_1,\ldots ,\rho_n),\alpha_1,\ldots ,\alpha_n)
\in \spec(n)$.
We say that the transposition $s_i$ is {\em admissible} for $\alpha$ if one
of the following conditions holds: 
 
\noi (i) $\rho_i\not=\rho_{i+1}$ or 

\noi (ii) $\rho_i=\rho_{i+1}$ and $\alpha_i \not=
\alpha_{i+1} \pm \frac{|G|}{\dim(V^{\rho_i})}$.

The following two observations are easy to see: 

\noi (a) For 
$\alpha, \beta \in \spec(n)$, we have $\alpha \sim \beta$ if $\alpha$ 
is obtained from $\beta$ by a 
sequence of admissible transpositions.

\noi (b) We have
\beq \nonumber
\lefteqn{((\rho_1,\ldots ,\rho_n),\alpha_1,\ldots ,\alpha_n) \in \spec(n)
\mbox{ implies }}\\ \label{rt} 
&&((\rho_1,\ldots ,\rho_{n-1}),\alpha_1,\ldots ,\alpha_{n-1}) \in \spec(n-1).
\eeq

\section{\large{Content vectors and Young $G$-tableaux}}

In the
Vershik-Okounkov theory Young $G$-tableaux are related to irreducible 
representations of $G_n$ via their content vectors. 
Let us define these first.

Let $\alpha =(a_1,a_2,\ldots ,a_n)\in {\mathbb Z}^n$. We say that
$\alpha$ is a {\em content vector} if

\noi (i) $a_1 = 0$.

\noi (ii) $\{a_i - 1 , a_i + 1\}\cap \{a_1,a_2,\ldots ,a_{i-1}\} \not=
\emptyset,\mbox{ for all }i> 1$.

\noi (iii) if $a_i = a_j = a$ for some $i<j$ then $\{a-1,a+1\}\subseteq
\{a_{i+1},\ldots ,a_{j-1}\}$ (i.e., between two occurrences of $a$ there
should also be occurrences of $a-1$ and $a+1$).

Condition (ii) in the definition above can be replaced (in the presence of
conditions (i) and (iii)) by condition (ii') below.

\noi (ii') 
For all $i > 1$, if $a_i > 0$ then $a_j = a_i - 1$ for some $j< i$ and if
$a_i < 0$ then $a_j = a_i + 1$ for some $j < i$.

The set of all content vectors of length $n$ is denoted $\cont(n) \subseteq
{\mathbb Z}^n$. It is convenient to assume that the empty sequence is a
content vector of length 0 and is the unique element of $\cont(0)$. 
 
Let $\alpha=((\rho_1,\ldots ,\rho_n),a_1,\ldots ,a_n)$, where $\rho_i\in
G^{\wedge}$ for all $i$, and $(a_1,\ldots ,a_n)\in \C^n$. 
For $\sigma\in G^{\wedge}$, let $\sigma(J)=\{j_1 < j_2 < \cdots <
j_{n_{\sigma,\alpha}}\} \seq \{1,2,\ldots ,n\}$ be the set of indices
satisfying $\rho_{j_i}=\sigma$, $i=1,\ldots ,n_{\sigma,\alpha}$ and
$\rho_{l}\not=\sigma$ for $l\in \{1,2,\ldots ,n\}-\sigma(J)$. 
Let $\sigma(\alpha)$ be the sequence 
$$\left( \frac{\dim(V^{\sigma})}{|G|}a_{j_1},\ldots 
, \frac{\dim(V^{\sigma})}{|G|}a_{j_{n_{\sigma,\alpha}}}\right).$$

We say that
$\alpha$ is a {\em content vector with respect to $G$ of length $n$} if   
$\sigma(\alpha)\in \cont(n_{\sigma,\alpha})$
for all $\sigma \in G^{\wedge}$. Since $\dim(V^{\sigma})$ divides $|G|$ it
follows that $a_i\in \Zi$ for all $i$.
Denote by $\contg(n)\seq \Zi^n$ the set of all content vectors with respect to $G$ of
length $n$.

\bt We have $\spec(n)\seq \contg(n)$.
\et
\pf 
Let $\alpha = (\rho, a_1,\ldots ,a_n)\in \spec(n)$, where
$\rho=(\sigma,\ldots ,\sigma)$, $\sigma \in G^{\wedge}$.
We will show that
$$\left( \frac{\dim(V^{\sigma})}{|G|}a_1,\ldots 
, \frac{\dim(V^{\sigma})}{|G|}a_n \right)\in
\cont(n).$$
Using 
Theorem \ref{mtl}(v) and (\ref{rt}) we see that this proves the result. 

Clearly $a_1 = 0$ as $X_1
= 0$. We verify conditions (ii) and (iii) in the definition of content
vectors by induction on $n$. 
Since $X_2 = b_1$ we have, from (\ref{bve}), that
$\frac{\dim(V^{\sigma})}{|G|}\;a_2 = \pm 1$ 
and thus condition (ii) is verified (and
condition (iii) does not apply). Now assume $n\geq 3$.

We first verify condition (ii). If $a_{n-1} = a_n \pm
\frac{|G|}{\dim(V^{\sigma})}$ 
there is nothing to prove, so assume this does not
hold.
Then the transposition $(n-1, n)$ is admissible for $\alpha$ and, by Theorem
\ref{mtl}(vi), 
$(\rho , a_1,\ldots ,a_{n-2},a_n,a_{n-1})\in \spec (n)$. 
Now, by (\ref{rt}) and the induction hypothesis, 
$\{a_n - \frac{|G|}{\dim(V^{\sigma})}, a_n + \frac{|G|}{\dim(V^{\sigma})}
\} \cap \{a_1,\ldots ,a_{n-2}\} \not= \emptyset$.
Thus condition (ii) is verified.

We now verify condition (iii). Now assume that $a_i=a_n=a$ for some $i<n$. We may assume that $i$ is the
largest possible index, i.e., $a$ does not occur between $a_i$ and $a_n$, so
$a\not\in \{a_{i+1},\ldots ,a_{n-1}\}$. Now assume that 
$a-\frac{|G|}{\dim(V^{\sigma})} \not\in
\{a_{i+1},\ldots ,a_{n-1}\}$. We shall derive a contradiction (the case
where $a+\frac{|G|}{\dim(V^{\sigma})} 
\not\in \{a_{i+1},\ldots ,a_{n-1}\}$ is similar).

By induction hypothesis the number $a+\frac{|G|}{\dim(V^{\sigma})}$ 
occurs in $\{a_{i+1},\ldots
,a_{n-1}\}$ at most once (for, if it occured twice, then  by the induction
hypothesis $a$ would also occur contradicting our choice of $i$). Thus there 
are two possibilities:
$$(a_i,\ldots ,a_n)=(a,*,\ldots ,*,a) \mbox{ or }
(a_i,\ldots ,a_n)
=(a,*,\ldots ,*,a+\frac{|G|}{\dim(V^{\sigma})},*,\ldots ,*,a),$$
where $*$ stands for a number different from $a-\frac{|G|}{\dim(V^{\sigma})}, 
a, a+\frac{|G|}{\dim(V^{\sigma})}$.

In the first case we can apply a sequence of admissible transpositions to
infer that $(\rho, \ldots, a,a, \ldots )\in \spec (n)$, contradicting Theorem  
\ref{mtl}(iii) and 
in the second case we can apply a sequence of admissible transpositions to
infer that $(\rho, \ldots, a,a+\frac{|G|}{\dim(V^{\sigma})}, a, \ldots )\in \spec (n)$, contradicting Theorem
\ref{mtl}(iv)(b).  $\Box$

Let $\alpha = ((\rho_1,\ldots ,\rho_n), a_1,\ldots ,a_n)\in \contg(n)$ 
We say that the transposition $s_i$ is {\em admissible} for $\alpha$
if $\rho_i\not=\rho_{i+1}$, or $\rho_i=\rho_{i+1}$ and $a_i \not = a_{i+1}
\pm \frac{|G|}{\dim(V^{\rho_i})}$. We   
define the following equivalence relation on $\contg(n)$: $\alpha \approx  
\beta$ if $\beta$ can be obtained from $\alpha$ by a sequence of (zero or  
more) admissible transpositions.

We now introduce Young $G$-tableaux into the picture.

Let $T_1\in \tabx(n)$ and assume that either $i$ and $i+1$ do not appear in
the same Young diagram of $T_1$ or that they are in the same Young diagram of
$T_1$ but do not appear in the
same row or same column of this Young diagram. 
Then exchanging $i$ and $i+1$ in $T_1$ produces
another
standard Young $G$-tableau $T_2 \in \tabx(n)$. We say that $T_2$ is obtained
from  
$T_1$ by an {\em admissible transposition}. For $T_1,T_2 \in \tabx(n)$,
define $T_1
\approx T_2$ if $T_2$ can be obtained from $T_1$ by a sequence of (zero or
more) admissible transpositions (it is easily seen that $\approx$ is an   
equivalence relation).

\bl \label{sl}
Let $T_1,T_2\in \tabx(n)$. Then $T_1 \approx T_2$ if and only 
$T_1$ and $T_2$ have the same shape.
\el
\pf
The only if part is obvious. To prove the if part we proceed as follows.
Let $\mu\in \Y_n(G^{\wedge})$. Enumerate the elements of $G^{\wedge}$ as
$\sigma_1,\ldots ,\sigma_t$. Let $n_i$ be the number of boxes in
the Young diagram $\mu(\sigma_i)$. Then $n_1+\cdots +n_t = n$.
 
Define the following element
$R$ of $\tabx(n,\mu)$: fill the Young diagram of $\mu(\sigma_1)$ with the
numbers $1,2,\ldots ,n_1$ in {\em row major order}, i.e., the first row with
the numbers $1,2,\ldots ,l_1$ (in increasing order, here $l_1 =$ 
length of first row), the second row
with $l_1+1,\ldots ,l_1+l_2$ (in increasing order, here 
$l_2 =$ length of second row) and so on till the last row of
$\mu(\sigma_1)$. Now fill the Young diagram of $\mu(\sigma_2)$ with the
numbers $n_1+1,\ldots ,n_1+n_2$ in row major order and so on till the last
Young diagram $\mu(\sigma_t)$.

We show that any $T \in \tabx(n,\mu)$ satisfies $T
\approx R$. This will prove the if part. Consider the last box of the last 
row of the last Young diagram $\mu(\sigma_t)$. 
Let $i$ be written in this box of $T$. Exchange $i$ and $i+1$
in $T$ (which is clearly an admissible transposition). Now repeat this
procedure with $i+1$ and $i+2$, then $i+2$ and $i+3$, and finally $n-1$ and
$n$. At the end of this sequence of admissible transpositions we have the  
number $n$ written in the last box of the last row of $\mu(\sigma_t)$. 
Now repeat the same procedure for $n-1,n-2,\ldots ,2$.  $\Box$

Let us make a remark about the proof of Lemma \ref{sl}. Let $s$ denote the
permutation that maps $R$ to $T$. 
Then the proof shows that $R$ can be
obtained from $T$ by a sequence of $\ell(s)$ admissible transpositions.
Thus $T$ can be
obtained from $R$ by a sequence of $\ell(s)$ admissible transpositions.

The {\em content}  
$c(b)$ of a box $b$ of a Young diagram is its 
$y$-coordinate $-$ its $x$-coordinate (our convention for drawing Young
diagrams is akin to writing down matrices with $x$-axis running downwards 
and $y$ axis running to the right). 
 
\bl \label{cl}
Let $\Phi : \tabx(n) \rar \contg(n)$ be defined as follows.
Given $T \in \tabx(n)$ and $1\leq i \leq n$, let $b_T(i)$ be the box (in one
of the Young diagrams of $T$) where the number $i$ resides.
Define
$$\Phi(T) = ((r_T(1),\ldots ,r_T(n)),\; 
\frac{|G|}{\dim(V^{r_T(1)})}\;c(b_T(1)),\ldots ,\;
\frac{|G|}{\dim(V^{r_T(n)})}\;c(b_T(n))).$$
Then $\Phi$ is a bijection which takes $\approx$-equivalent standard Young
$G$-tableaux to $\approx$-equivalent content vectors with respect to $G$.
\el
\pf The general case clearly follows from the $|G|=1$ case for which we
need to give a bijection between content vectors of 
length $n$ and standard Young tableaux with $n$ boxes. This is well known. 
The content vector of any standard Young tableau clearly satisfies
conditions (i), (ii), and (iii) in the definition of a content vector and
these conditions uniquely determine the numbers to be filled in the boxes of
the Young diagram. This bijection clearly preserves the $\approx$ relation.
$\Box$

\bt \label{ts}
(i) $\spec(n) = \contg(n)$ and the equivalence relations $\sim$ and
$\approx$ coincide.

\noi (ii) The map $\Phi^{-1}: \spec(n) \rar \tabx(n)$ is a
bijection and, for $\alpha ,\beta \in \spec(n)$, we have $\alpha \sim \beta$
if and only if $\Phi^{-1}(\alpha)$ and $\Phi^{-1}(\beta)$ have the same
shape.

\et
\pf We have

\noi (a) $\spec(n) \subseteq \contg(n)$.

\noi (b) If $\alpha \in \spec(n)$, $\beta \in \contg(n)$, and $\alpha \approx
\beta$ then it is easily seen that $\beta \in \spec(n)$ and $\alpha \sim  
\beta$. It follows that given an $\sim$-equivalence class $\cal A$ of $\spec
(n)$ and an $\approx$-equivalence class $\cal B$ of $\contg(n)$, either
$\cal A \cap \cal B = \emptyset$ or $\cal B \subseteq \cal A$.

\noi (c) $ |(\spec(n) / \sim)| = |G_n^{\wedge}| = |\pc_n(G_*)|=
|\Y_n(G^{\wedge})|$, 
since the number of irreducible
$G_n$-representations is equal to the number of conjugacy classes in $G_n$
and similarly for $G$.

\noi (d) $|(\contg(n) / \approx)| = |\Y_n(G^{\wedge})|$, by Lemmas \ref{cl} and \ref{sl}.

The four statements above imply part (i). Part (ii) is now
clear. $\Box$

Using Theorem \ref{ts} we may parametrize the irreducible representations of
$G_n$ by elements of $\Y_n(G^{\wedge})$. 
The following result is a reformulation of the GZ-decomposition in terms of
standard Young $G$-tableaux.
\bt \label{gzt}
Let $\mu\in \Y_n(G^{\wedge})$. Then we may index the GZ-subspaces of
$V^{\mu}$ by standard Young $G$-tableaux of shape $\mu$ and write the GZ-decomposition
(\ref{gzdi}) as
\beq \label{gzd} 
V^{\mu} &=& \oplus_{T\in \tabx(n,\mu)} V_T,
\eeq
where each $V_T$ is closed under the action of $G^n=G\times \cdots \times G$
($n$ factors) and, as a $G^n$-module, is isomorphic to the irreducible
$G^n$-module
$$V^{r_T(1)}\otimes V^{r_T(2)}\otimes \cdots \otimes V^{r_T(n)}.$$  
For $i=1,\ldots ,n$, the eigenvalue of $X_i$ on $V_T$ is given by
$\frac{|G|}{\dim(V^{r_T(i)})}\;c(b_T(i))$. $\Box$ 
\et

The branching rule for the pair $G_{n-1}\seq G_n$ is now clear. 
\bt \label{br}
Let $\mu\in \Y_{n+1}(G^{\wedge})$. Then we have a $G_n$-module isomorphism
\beqn
V^{\mu} &\cong& \oplus_{\sigma\in G^{\wedge}} \dim(V^{\sigma})\;\left( 
\oplus_{\lambda \in \mu \da \sigma} V^{\lambda} \right).\;\;\Box
\eeqn
\et
The dimension of an irreducible representation of $G_n$ easily
follows from Theorem
\ref{gzt}. For a Young diagram $\mu$ let $f^{\mu}$ denote the number of
standard Young tableaux of shape $\mu$.

\bt \label{dt}
Let $\mu\in \Y_n(G^{\wedge})$. Write the elements of $G^{\wedge}$ as
$\{\sigma_1,\ldots ,\sigma_t\}$ and set 
$$\mu_i = \mu(\sigma_i),\;n_i =
|\mu_i|,\;d_i=\dim(V^{\sigma_i}),\;i=1,\ldots ,t.$$
Then
$$\dim(V^{\mu}) = \binom{n}{n_1,\ldots ,n_t}\;f^{\mu_1}\cdots
f^{\mu_t}\;d_1^{n_1}\cdots d_t^{n_t}. \Box$$
\et

We now discuss the choice of a basis of $V^{\mu},\;\mu\in \Y_n(G^{\wedge})$,
with respect to which we may write down the matrices for the action of the
Coxeter generators $s_1,\ldots ,s_{n-1}$. We begin with an observation.

Fix $\mu\in \Y_n(G^{\wedge})$ and consider the irreducible $G_n$-module
$V^{\mu}$. Let $T\in \tabx(n,\mu)$ and 
let $p_T$ denote the projection of $V^{\mu}$ onto $V_T$ determined by the
decomposition (\ref{gzd}). 
Let
$s_i$ be an admissible transposition for $T$. Two cases arise:

\noi (a) $i$ and $i+1$ are in different Young diagrams of $T$: It follows
from Theorem \ref{mtl} (v) that 
$p_{s_i\cdot T}(s_i\cdot B)$ is a basis of $V_{s_i\cdot
T}$ for any basis $B$ of $V_T$.

\noi (b) $i$ and $i+1$ are in the same Young diagram of $T$ but are not in
the same row or the same column of this Young diagram: 
Let $0\not= v\in V_T$. It follows from
(the proof of) Theorem \ref{mtl} (vi) that
$s_i \cdot v$ is the  sum of a nonzero rational multiple of
$\tau_{i,V_T}(v)$ and a nonzero vector in $V_{s_i \cdot T}$ and that  
the map $V_T \rar V_{s_i\cdot T}$ given by $v\mapsto p_{s_i\cdot T}(s_i\cdot
v)$ is  a linear isomorphism. In particular,
$p_{s_i\cdot T}(s_i\cdot B)$ is a basis of $V_{s_i\cdot
T}$ for any basis $B$ of $V_T$.  

Now let $R$ be the tableau defined in the
proof of Lemma \ref{sl}. 
Fix a basis $B_R$ of $V_R$. 
Consider
a standard $G$-tableau $T \in \tabx(n,\mu)$. 
Let $s$ be the permutation that maps $R$ to
$T$. 
Define
\beq \label{bd}
B_T &=& \{ p_T(s\cdot v)\;|\; v\in B_R\},
\eeq
and define $\ell(T)$, the {\em length} of $T$, to be  $\ell(s)$.
The following result now easily follows, by induction on the length of $T\in
\tabx(n,\mu)$, using observations (a), (b) above and the fact, remarked
after the proof of Lemma \ref{sl}, that $T\in \tabx(n,\mu)$ can be obtained
from $R$ by a sequence of $\ell(T)$ admissible transpositions and no fewer
Coxeter transpositions.

\bl \label{am}
(i) $B_T$ is a basis of $V_T$, for all $T\in \tabx(n,\mu)$.

\noi (ii) Let $T\in \tabx(n,\mu)$ and let $s_i$ be an admissible transposition
for $T$. Then

(a) If $i$ and $i+1$ are in different Young diagrams of $T$ we have 
$$B_{s_i\cdot T} = \{ p_{s_i\cdot T}(s_i\cdot v)\;|\;v\in B_T\}.$$

(b) If $i$ and $i+1$ are in the same Young diagram of $T$ but not in the
same row or same column of this Young diagram we have
\beqn
B_{s_i\cdot T} &=& \{ p_{s_i\cdot T}(s_i\cdot v)\;|\;v\in B_T\},\;\;\mbox{if
} \ell(s_i\cdot T)=\ell(T)+1,\\
B_{s_i\cdot T} &=& \{(1-r^{-2}) p_{s_i\cdot T}(s_i\cdot v)\;|\;v\in B_T\},\;\;\mbox{if
} \ell(s_i\cdot T)=\ell(T)-1,
\eeqn
where $r=\frac{(c(b_T(i+1))-c(b_T(i)))}{|G|}\;\mbox{dim}(V^{r_T(i)}).$
\el
We now choose a basis of $V_R$ in a certain way and then apply the method
above to get bases of all the GZ-subspaces. 
For $\sigma \in G^{\wedge}$, fix a
basis $B^{\sigma}$ of $V^{\sigma}$. Then, for $\rho=(\rho_1,\ldots ,\rho_n)$,
where $\rho_i \in G^{\wedge}$ for all $i$, we have that
$B^{\rho}=B^{\rho_1}\otimes \cdots \otimes B^{\rho_n}$ is a basis of
$V^{\rho}=V^{\rho_1}\otimes \cdots \otimes V^{\rho_n}$. 
Thus, for $T\in
\tabx(n,\mu)$, we have that $B^{r_T}$ is a basis of $V^{r_T}$, where
$r_T=(r_T(1),\ldots ,r_T(n))$. 
Let $N_{i,r_T}$ be
the matrix of the switch operator $\tau_{i,r_T}$ on $V^{r_T}$ with respect
to the basis $B^{r_T}$.

Let $R$ be the standard $G$-tableau defined above and fix a $G$-linear
isomorphism 
$$f: V^{r_R(1)}\otimes \cdots \otimes V^{r_R(n)}\rar V_R.$$
Define the basis
$$\ol{B_R} = f(B^{r_R(1)}\otimes \cdots \otimes B^{r_R(n)})$$
of $V_R$. 
Now use (\ref{bd}) to define a basis $\ol{B_T}$ of $V_T$ for all
$T\in \tabx(n,\mu)$.

Let $T\in \tabx(n,\mu)$ and $s$ be the permutation that maps $R$ to $T$. Now
$S_n$ acts on $V^{r_T}$ by permuting the coordinates and the image of the
action of $s^{-1}$ on $V^{r_T}$ is $V^{r_R}$. The following result now
follows. 
\bl \label{am1}
The map $V^{r_T} \rar V_T$ given by $v\mapsto p_T(sfs^{-1}(v))$ is a
$G^n$-linear isomorphism that takes the basis $B^{r_T}$ of $V^{r_T}$ to the
basis $\ol{B_T}$ of $V_T$.
Thus the
matrix $N_{i,T}$ of $\tau_{i,V_T}$ with respect to $\ol{B_T}$ is equal to
$N_{i,r_T}$.
\el

We now have the following result.
\bt Consider the basis $\cup_{T\in \tabx(n,\mu)} \ol{B_T}$  
of $V^{\mu}$.
Fix $T\in \tabx(n,\mu)$ and let 
$\Phi(T) = ((\rho_1,\ldots ,\rho_n), a_1,\ldots , a_n)$. 
Let $s_i$
be a Coxeter generator. 
Let $I$ denote the $|\ol{B_T}|\times |\ol{B_T}|$ identity matrix.  
Set $r=\frac{(a_{i+1}-a_i)\mbox{dim}(V^{\rho_i})}{|G|}$ and $N=N_{i,T}$. 

The action of $s_i$ on $V_T$ is as follows.

\noi (i) If $i$ and $i+1$ are in the same column of the same Young diagram
of $T$ 
then $V_T$ is closed under the action of $s_i$ and the matrix of this action
with respect to the basis $\ol{B_T}$ is $N$.

\noi (ii) If $i$ and $i+1$ are in the same row of the same Young diagram
of $T$ 
then $V_T$ is closed under the action of $s_i$ and the matrix of this action
with respect to the basis $\ol{B_T}$ is $-N$.

\noi (iii) Suppose $i$ and $i+1$ are not in the same Young diagram of $T$. 
Let $S=s_i\cdot T$. 
Then $V_T \oplus V_S$ is closed under the action of
$s_i$ and the  matrix of this 
action, with respect to the basis $\ol{B_T} \cup \ol{B_S}$, is given by
$$\left[ \ba{cc}  0 & I \\I & 0 \ea \right].$$

\noi (iv) Suppose $i$ and $i+1$ are in the same Young diagram of $T$ but 
not in the same row or same column of this Young diagram. Let
$S = s_i \cdot T$. Then $N=N_{i,S}$.

If $\ell(S) = \ell (T) + 1$ then 
$V_T \oplus V_S$ is closed under the action of
$s_i$ and the  matrix of this 
action, with respect to the basis $\ol{B_T} \cup \ol{B_S}$, is given by
$$  
\left[\ba{cc}
r^{-1}N & (1 - r^{-2}) I \\ 
I & - r^{-1}N \ea \right].$$

If $\ell(S) = \ell (T) - 1$ then the matrix of 
the action of $s_i$ on the subspace $V_T \oplus V_S$
with respect to the basis basis $\ol{B_T} \cup \ol{B_S}$ is given by 
the transpose of the matrix above.
\et
\pf Parts (i), (ii), (iii) and part (iv) with $\ell(S)=\ell(T)+1$ follow
from Theorem \ref{mtl}, Lemma \ref{am}, and Lemma \ref{am1} above.
To prove the case $\ell(S)=\ell(T)-1$ of part (iv), switch $T$ and $S$ in
the $\ell(S)=\ell(T)+1$ case along with switching $a_i$ and $a_{i+1}$. This
is equivalent to transposing the matrix. $\Box$

The basis of $V^{\mu}$ and the action of $s_i$ 
described above correspond to Young's {\em seminormal
form} in the case of the symmetric groups.
Now let us consider the analog of Young's {\em orthogonal form}.
Since $V^{\mu}$ is irreducible there is
a unique (upto scalars) $G_n$-invariant inner product on $V^{\mu}$. Choose
and fix one such inner product. Since
the branching from $H_{i,n}$ to $H_{i-1,n}$ is multiplicity free we have
that
the decomposition of an irreducible $H_{i,n}$-module into irreducibles
$H_{i-1,n}$-modules is orthogonal. It follows that the GZ-decomposition
(\ref{gzd}) of $V^{\mu}$ is orthogonal. 

For $\sigma \in G^{\wedge}$, fix a $G$-invariant inner product (unique upto
scalars) on $V^{\sigma}$, and fix an orthonormal
basis $C^{\sigma}$ of $V^{\sigma}$. Then, for $\rho=(\rho_1,\ldots ,\rho_n)$,
where $\rho_i \in G^{\wedge}$ for all $i$, we have that
$C^{\rho}=C^{\rho_1}\otimes \cdots \otimes C^{\rho_n}$ is an orthonormal 
basis of
$V^{\rho}=V^{\rho_1}\otimes \cdots \otimes V^{\rho_n}$ (under the inner
product obtained by multiplying the component inner products). 
Thus, for $T\in
\tabx(n,\mu)$, we have that $C^{r_T}$ is an orthonormal 
basis of $V^{r_T}$, where
$r_T=(r_T(1),\ldots ,r_T(n))$. 
Let $M_{i,r_T}$ be
the matrix of the switch operator $\tau_{i,r_T}$ on $V^{r_T}$ with respect
to the basis $C^{r_T}$.

Let $R$ be the standard $G$-tableau defined above and fix a $G^n$-linear
isometry 
$$f: V^{r_R(1)}\otimes \cdots \otimes V^{r_R(n)}\rar V_R.$$
Define the orthonormal basis
$$C_R = f(C^{r_R(1)}\otimes \cdots \otimes C^{r_R(n)})$$
of $V_R$. 
Now use (\ref{bd}) to define a basis $C_T$ of $V_T$ for all
$T\in \tabx(n,\mu)$.

Let $T\in \tabx(n,\mu)$ and $s$ be the permutation that maps $R$ to $T$. 
We have 
\bl \label{am2}
The map $V^{r_T} \rar V_T$ given by $v\mapsto p_T(sfs^{-1}(v))$ is a
$G^n$-linear isometry that takes the basis $C^{r_T}$ of $V^{r_T}$ to the
basis $C_T$ of $V_T$.
Thus the
matrix $M_{i,T}$ of $\tau_{i,V_T}$ with respect to $C_T$ is equal to
$M_{i,r_T}$.
\el

The following result can be proved along the lines of the 
previous result.
\bt Consider the orthonormal basis $\cup_{T\in \tabx(n,\mu)} C_T$  
of $V^{\mu}$ defined above. 
Fix $T\in \tabx(n,\mu)$ and let 
$\Phi(T) = ((\rho_1,\ldots ,\rho_n), a_1,\ldots , a_n)$. 
Let $s_i$
be a Coxeter generator. 
Let $I$ denote the $|C_T|\times |C_T|$ identity matrix and let 
$M_{i,T}$ denote the matrix of $\tau_{i,V_T}$ with respect to the basis $C_T$. 
Set $r=\frac{(a_{i+1}-a_i)\mbox{dim}(V^{\rho_i})}{|G|}$ and $M=M_{i,T}$. 

The action of $s_i$ on $V_T$ is as follows.

\noi (i) If $i$ and $i+1$ are in the same column of the same Young diagram
of $T$ 
then $V_T$ is closed under the action of $s_i$ and the matrix of this action
with respect to the basis $C_T$ is $M$.

\noi (ii) If $i$ and $i+1$ are in the same row of the same Young diagram
of $T$ 
then $V_T$ is closed under the action of $s_i$ and the matrix of this action
with respect to the basis $C_T$ is $-M$.

\noi (iii) Suppose $i$ and $i+1$ are not in the same Young diagram of $T$. 
Let $S=s_i\cdot T$. 
Then $V_T \oplus V_S$ is closed under the action of
$s_i$ and the  matrix of this 
action, with respect to the basis $C_T \cup C_S$, is given by
$$\left[ \ba{cc}  0 & I \\I & 0 \ea \right].$$

\noi (iv) Suppose $i$ and $i+1$ are in the same Young diagram of $T$ but 
not in the same row or same column of this Young diagram. Let
$S = s_i \cdot T$. Then $M=M_{i,S}$. 

Then $V_T \oplus V_S$ is closed under the action of
$s_i$ and the  matrix of this 
action, with respect to the basis $C_T \cup C_S$, is given by
$$  
\left[\ba{cc}
r^{-1}M & \sqrt{1 - r^{-2}}\;I \\ 
\sqrt{1-r^{-2}}\;I & - r^{-1}M \ea \right].$$
\et

\section{\large{Generalized Johnson scheme}}

The simplest nontrivial 
examples of the Vershik-Okounkov theory are the classical ``Johnson schemes" 
and the ``generalized Johnson
schemes" of Ceccherini-Silberstein, Scarabotti, and Tolli
{\bf\cite{cst1,cst2}} (also see {\bf\cite{ms}}). We consider multiplicity
free $S_n$, $G_n$-actions and explicitly write down the GZ-vectors (in the
$S_n$ case) and the GZ-subspaces (in the $G_n$ case) and also identify the
irreducibles which occur. 

We begin with the $S_n$ action.
Let $B(n)$ denote the set of all subsets of 
$[n]= \{1,2,\ldots ,n\}$ and, for $0\leq i \leq n$, let $B(n)_i$ denote the
set of all subsets of $[n]$ with cardinality $i$. 
There is a natural action
of $S_n$ on $B(n)_i$ and $B(n)$.
For a finite set $S$, let $V(S)$ denote the complex vector space with $S$ as
basis.

We have the
following direct sum decomposition into $S_n$-submodules of the permutation
representation of $S_n$ on $V(B(n))$:
\beq 
V(B(n)) &=& V(B(n)_0) \oplus V(B(n)_1) \oplus \cdots \oplus
V(B(n)_n).
\eeq

The following result is classical ({\bf\cite{cst2,st1}}). We give a
constructive proof that produces an explicit Gelfand-Tsetlin basis.

\bt \label{js}
For $0\leq i \leq n$, $V(B(n)_i)$ is a multiplicity free $S_n$-module
with $S_n$-module isomorphism
\beqn
V(B(n)_i) &\cong & \bigoplus_k \; V^{(n-k,k)},
\eeqn
where the sum is over all partitions $(n-k,k)$ of $n$ with atmost two parts
satisfying $k\leq i \leq n-k$.

\et
An element $v\in V(B(n))$ is {\em homogeneous} if $v\in V(B(n)_k)$ for some
$k$. We say that a nonzero homogeneous element $v$ is of {\em rank}
$k$, and we write $r(v)=k$, if $v\in V(B(n)_k)$. 
The {\em
up operator}  $U_n:V(B(n))\rar V(B(n))$ is defined, for $X\in B(n)$, by
$$U_n(X)= \sum_{Y} Y,$$
where the sum is over all $Y\in B(n)$ covering $X$, i.e., $X\subseteq Y$ and
$|Y|=|X|+1$.

A {\em symmetric Jordan chain} (SJC) in $V(B(n))$ is a sequence   
$v=(v_1,\ldots ,v_h)$ of nonzero homogeneous elements of $V(B(n))$
such that $U_n(v_{i-1})=v_i$, for    
$i=2,\ldots h$, $U_n(v_h)=0$, and    
$r(v_1) + r(v_h) = n$, if $h\geq   
2$, or else $2r(v_1)= n$, if $h=1$.
Note that the
elements of the sequence $v$ are linearly independent, being nonzero and of
different ranks. We say that $v$ {\em starts} at rank $r(v_1)$ and {\em
ends} at rank $n-r(v_1)$. We do not distinguish between the sequence
$(v_1,\ldots ,v_h)$ and the underlying set $\{v_1,\ldots ,v_h\}$.
A {\em symmetric Jordan basis} (SJB)  of $V(B(n))$ is a basis of $V(B(n))$
consisting of a disjoint union of SJC's in $V(B(n))$. Given an SJB $J(n)$ of
$V(B(n))$ and $0\leq k \leq n/2$, let $\J(n,k)$ denote the set of all SJC's
in $J(n)$ starting at rank $k$ and ending at rank $n-k$ and let $J(n,k)$
denote the union of all SJC's in $\J(n,k)$.

Given $T\in \tab(n,\mu)$, where $\mu$ has atmost two rows, we denote by $T
+_1 (n+1)$ the standard Young tableaux obtained from $T$ by adding $n+1$ at
the end of the first row.
Similarly, given $T\in \tab(n,\mu)$, where $\mu$ has atmost two rows with
the second row containing fewer elements than the first row, we denote by $T
+_2 (n+1)$ the standard Young tableaux obtained from $T$ by adding $n+1$ at
the end of the second row.
The basic idea of the following algorithm is from {\bf\cite{sr1}}, though we
have added new elements here, namely, Theorems \ref{cgzb1} and \ref{ev}. 
\bt \label{cgzb}
There exists an inductive procedure to explicitly construct an SJB $J(n)$ of
$V(B(n))$ and, for $0\leq k \leq n/2$, 
a bijection 
\beq \label{tsjc}
&B_{n,k} : \tab(n,(n-k,k)) \rar \J(n,k).&
\eeq
\et
\pf The case $n=1$ is clear.

Let $V=V(B(n+1))$. Define $V(0)$ to be the subspace of $V$ generated by all
subsets of $[n+1]$ not containing $n+1$ and define $V(1)$ to be the subspace
of $V$ generated by all subsets of $[n+1]$ containing $n+1$. We have
$V=V(0)\oplus V(1)$. The linear map $\R:V(0)\rar
V(1)$, given by  $X\mapsto X\cup \{n+1\},\;X\subseteq [n]$ is an
isomorphism.
We write $\R(v) = \ol{v}$. We write $U$ for the up operator $U_{n+1}$
on $V$ and we write $U_0$ for the up
operator on $V(0)$ ($=V(B(n))$). We have, for $v\in V(0)$, 
\beq \label{a} 
&U(v)=U_0(v)+\ol{v},\;\;U(\ol{v})=\ol{U_0(v)}.
\eeq

By induction hypothesis there is an SJB $J(n)$ of $V(B(n))=V(0)$ and
bijections $B_{n,k}$ as in (\ref{tsjc}) above.
We shall now produce an SJB $J(n+1)$ of $V$ by producing, for each
SJC in $J(n)$, either one or two SJC's in $V$ such that the   
collection of all these SJC's is a basis. 

Let $0\leq k \leq n/2$. Consider $T\in \tab(n,(n-k,k))$ and
consider the SJC $B_{n,k}(T)=(x_k,\ldots ,x_{n-k})\in \J(n,k)$, 
where $r(x_k) = k$.    

We now consider two cases.

(a) $k=n-k$ : From (\ref{a}) we have $U(x_k)=\ol{x_k}$ and
$U(\ol{x_k})=\ol{U_0(x_k)}=0$. Since $\R$ is an isomorphism $\ol{x_k}\not=0$.
Define
\beq \label{ba2}
B_{n+1,k}(T +_1 (n+1)) &=&  (x_k, \ol{x_k}).
\eeq 

(b) $k< n-k$ : Set $x_{k-1}=x_{n+1-k}=0$ and define
\beq \label{ba3}
&y=(y_k,\ldots ,y_{n+1-k}), \mbox{ and }
z=(z_{k+1},\ldots ,z_{n-k}), \eeq
by
\beq \label{ba4}
y_l &=& x_l + (l-k)\, \ol{x_{l-1}},\;\;k\leq l \leq n+1-k. \\\label{ba5}
z_l &=& (n-k-l+1)\,\ol{x_{l-1}} -x_l,\;\;k+1\leq l \leq n-k.\eeq

From (\ref{a}) we have
\beq \label{b}
&U(\ol{x_l})=\ol{U_0(x_l)} = \ol{x_{l+1}},\;k\leq l\leq n-k &
\eeq 

It thus follows from (\ref{a}) and (\ref{b}) that, for $k\leq l < n+1-k$, we
have
$$U(y_l)=U(x_l + (l-k)\ol{x_{l-1}}) = x_{l+1} + \ol{x_l} +
(l-k)\ol{x_l}=x_{l+1}+(l-k+1)\ol{x_l}=y_{l+1}.$$ 
Note that when $l=k$ the second step above is justified
because of the presence of the $(l-k)$ factor even though
$U(\ol{x_{k-1}})=0\not= \ol{x_k}$. We
also have $U(y_{n+1-k})=U((n+1)\ol{x_{n-k}})=(n+1)\ol{U_0(x_{n-k})}=0$.

Similarly, for $k+1\leq l < n-k$, we have
$$U(z_l)=U((n-k-l+1)\ol{x_{l-1}}-x_l) = (n-k-l+1)\ol{x_l} -x_{l+1} -
\ol{x_l} = (n-k-l)\ol{x_l} - x_{l+1} = z_{l+1}.$$ 
and $U(z_{n-k})=U(\ol{x_{n-k-1}} - x_{n-k})=\ol{x_{n-k}}-\ol{x_{n-k}}=0.$

Since $y_k = x_k\not= 0$, $y_{n+1-k} = (n+1)\ol{x_{n-k}}\not= 0$,
$x_l$ and $\ol{x_{l-1}}$ are linearly independent, for $k+1\leq l \leq n-k$
and the $2\times 2$ matrix 
$$\left(\ba{rl} 1 & l-k \\ 
                           -1 & n-k-l+1 \ea\right)$$
is nonsingular for $k+1\leq l \leq n-k$, it follows that (\ref{ba3}) gives
two independent SJC's in $V$. 
Define 
\beqn
B_{n+1,k}(T +_1 (n+1)) &= &y,\\
B_{n+1.k+1}(T +_2 (n+1))&=&z,
\eeqn
and set $J(n+1)$ to be the union of all SJC's obtained in steps (\ref{ba2})
and (\ref{ba3}) above.

Since $V=V(0)\oplus V(1)$ and $\R$ is an isomorphism it follows that
$J(n+1)$ is an SJB of $V$. That the maps $B_{n+1,k}$
are bijections is also clear. $\Box$ 

\bex \label{example}
{\em
In this example we work out the SJB's of $V(B(n))$, for $n=2,3$, starting
with the SJB of $V(B(1))$, using the formulas 
(\ref{ba2}, \ref{ba3}, \ref{ba4}, \ref{ba5})
given in the proof of Theorem \ref{cgzb}.

\noi (i) The SJB of $V(B(1))$ is given by
\beqn
&(\,\,\emptyset \,,\,\{1\}\,\,)&
\eeqn

\noi (ii) The SJB of $V(B(2))$ consists of
\beqn
&(\,\,\emptyset \,,\, \{1\} + \{2\}\,,\, 2\{1,2\}\,\,)&\\
&(\,\,\{2\} - \{1\}\,\,)&
\eeqn

\noi (iii) The SJB of $V(B(3))$ consists of
\beqn
&(\,\,\emptyset \,,\, \{1\} + \{2\} + \{3\}\,,\,
2(\{1,2\}+\{1,3\}+\{2,3\})\,,\,6\{1,2,3\}\,\,)&\\
&(\,\,2\{3\}-\{1\}-\{2\}\,,\,\{1,3\}+\{2,3\}-2\{1,2\}\,\,)&\\
&(\,\,\{2\}-\{1\}\,,\,\{2,3\}-\{1,3\}\,\,)&
\eeqn

}\eex

For $0\leq k\leq i \leq n-k \leq n$ define
\beqn
J(n,k,i) &=& \{ v\in J(n,k) : r(v)=i \}.
\eeqn
Let $W(n,k,i)$ be the subspace of $V(B(n)_i)$ spanned by $J(n,k,i)$. 
Then we have the
direct sum decomposition
\beq \label{irred}
V(B(n)_i) &=& \bigoplus_{k=0}^{\mbox{min}\{i,n-i\}} W(n,k,i),\;\;
0\leq i \leq n.
\eeq
We claim that each $W(n,k,i)$ is a $S_{n}$-submodule of $V(B(n)_i)$. 
We prove this by
induction on $i$, the case $i=0$ being clear. Assume inductively that
$W(n,0,i-1),\ldots ,W(n,i-1,i-1)$ are submodules, where $i <  \lfloor n/2
\rfloor$.
Since $U_n$ is $S_{n}$-linear, $U_n(W(n,j,i-1))=W(n,j,i)$, 
$0\leq j \leq i-1$ are
submodules. Now consider $W(n,i,i)$.
Let $u\in W(n,i,i)$ and $\pi\in S_{n}$. Since $U_n$ is
$S_{n}$-linear we have $U_n^{n+1-2i}(\pi u)=\pi U_n^{n+1-2i}(u)=0$.
It follows that
$\pi u \in W(n,i,i)$. So the claim is proven for $0\leq i \leq n/2$ and it
follows for $i> n/2$ since $U_n$ is $S_n$-linear.

\bt \label{cgzb1}
As $S_n$-modules we have
\beq \label{irrep}
W(n,k,i) &\cong& V^{(n-k,k)},\;\;0\leq k\leq i \leq n-k \leq n.
\eeq
\et
\pf By induction on $n$. The cases $n=1,2,3$ can be directly verified from
Example \ref{example} (the main point to check is that $W(3,1,1)$ is the
standard representation of $S_3$).

Now assume we have proven the result upto $n\geq 3$.
By the algorithm of Theorem \ref{cgzb} we have, for $0\leq i \leq n+1$,
\beq \label{s0} 
W(n+1,k,i) &=& W(n,k,i) \oplus \ol{W}(n,k-1,i-1),
\eeq
where $\ol{W}=\{ \ol{v}\;|\; v\in W(n,k-1,i-1)\}$ (in the notation used in
the proof of Theorem \ref{cgzb}) and where $W(n,k,i)$ is
taken to be the zero subspace if $i< k$ or $i > n-k$.

Now, $W(n+1,k,i)$ is a $S_{n+1}$-module 
and it is easily seen that
$W(n,k,i)$ and $\ol{W}(n,k-1,i-1)$ are $S_n$-submodules of $W(n+1,k,i)$. By
induction hyphothesis we have, as $S_n$-modules,
\beq \label{s1}
W(n,k,i) &=& V^{(n-k,k)},\\ \label{s2}
\ol{W}(n,k-1,i-1) &=& V^{(n-k+1,k-1)}.
\eeq
Suppose an $S_{n+1}$-irreducible
$V^{\lambda}$, where the Young diagram $\lambda$ has 3 or more rows, occurs
in $W(n+1,k,i)$. Since $n+1\geq 4$, it follows that $\lambda$ has an inner
corner whose removal still leaves 3 or more rows. By the branching rule this
contradicts (\ref{s0}), (\ref{s1}), and (\ref{s2}). So, for any $S_{n+1}$-irreducible
$V^{\lambda}$ occuring in $W(n+1,k,i)$, there are atmost two rows in
$\lambda$.  It is now easy to see using the branching rule and (\ref{s1})
and (\ref{s2}) that $W(n+1,k,i) \cong V^{(n+1-k,k)}$. $\Box$

Theorem \ref{js} now follows from (\ref{irred}) and 
Theorem \ref{cgzb1}.
We also have that 
\beq \label{dimi}
\dim(V^{(n-k,k)})=\binom{n}{k} - \binom{n}{k-1}.
\eeq
Summing (\ref{irred}) over $i$ and taking dimensions we get 
\beq \label{bi}
2^n &=& \sum_{k=0}^{\lfloor n/2\rfloor} (n-2k+1)\;\left\{\binom{n}{k} -
\binom{n}{k-1}\right\}.
\eeq 

We denote the YJM elements of $S_n$ by $Y_1,\ldots ,Y_n$.

\bt \label{ev}
For $T\in \tab(n,(n-k,k))$ and every vector $v$ in the SJC 
$B_{n,k}(T)$ we have 
\beq \label{tsjc1} 
Y_j(v) &=& c(b_T(j))v,\;j=1,2,\ldots ,n.
\eeq
\et
\pf We first show inductively that each element
of $J(n)$ is a simultaneous eigenvector of $Y_1,\ldots ,Y_n$, the case $n=1$
being clear. 

Note that if $v\in V(B(n)_k)$ is an eigenvector for
$Y_i$, for some $1\leq  i \leq n$,
then $\ol{v}\in V(B(n+1)_{k+1})$ is also an eigenvector for $Y_i$ with the
same eigenvalue. Thus it follows from
(\ref{ba2}, \ref{ba3}, \ref{ba4}, \ref{ba5})
that each element of $J(n+1)$ is an eigenvector for $Y_1,\ldots ,Y_n$. It
remains to show that each element of $J(n+1)$ is an eigenvector for
$Y_{n+1}$.

We now have from Theorem \ref{cgzb1} that, for $0\leq i \leq \frac{n+1}{2}$,
$W(n+1,0,i),\ldots ,W(n+1,i,i)$ are mutually
nonisomorphic irreducibles. 
Consider the $S_{n+1}$-linear map
$f:V(B(n+1)_i)\rar
V(B(n+1)_i)$ given by $f(v) = av$, where
$$ a = \mbox{ sum of all transpositions in $S_{n+1}$ }= Y_1+\cdots
+Y_{n+1}.$$
It follows by Schur's lemma that there exist scalars $\alpha_0 ,\ldots
,\alpha_i$ such that $f(u)=\alpha_k u$, for $u\in W(n+1,k,i)$. Thus each
element
of $J(n+1,k,i)$ is an eigenvector for $Y_1+\cdots +Y_{n+1}$ (and also for
$Y_1,\ldots ,Y_n$). It follows that each element of $J(n+1,k,i)$ is an   
eigenvector for $Y_{n+1}$. 

The paragraph above has shown that the first element of each symmetric
Jordan chain in $J(n+1)$ 
is a simultaneous  eigenvector for $Y_1, \ldots ,Y_{n+1}$.
It now follows (since $U_{n+1}$ is $S_{n+1}$-linear) 
that each element of $J(n+1)$
is a simultaneous  eigenvector for $Y_1, \ldots ,Y_{n+1}$.

We are left to show that, for each $v\in B_{n,k}(T)$, the eigenvalues of
$Y_1,\ldots ,Y_n$ on $v$ are given by (\ref{tsjc1}). We can show this by
induction, the case $n=1$ being trivial.

Just like above the eigenvalues of $Y_1,\ldots ,Y_n$ on $v\in B_{n+1,k}(T)$ 
will continue to satisfy (\ref{tsjc1}). Now, since $v$ is an eigenvector for
$Y_{n+1}$ and $v$ lies in an $S_{n+1}$-irreducible isomorphic to
$V^{(n+1-k,k)}$, it follows that the eigenvalue of $Y_{n+1}$ on $v$ also
satisfies (\ref{tsjc1}).
That completes the proof. $\Box$

See {\bf\cite{f}} for an elegant direct construction of the GZ-basis given above
and see {\bf\cite{fkmw}} for an application to complexity theory.

Now we study the $G_n$ analog of the $S_n$-action considered above,
defined in {\bf\cite{cst1,cst2}}.
Let $G$ be a finite group acting on the finite set $X$. Assume that the
corresponding permutation representation on $V(X)$ is multiplicity free.
This implies, in
particular, that the action is transitive.

Let $L_0$ be a symbol not in $X$ and let $Y$ denote the {\em alphabet}
$Y=\{L_0\}\cup X$. We call the elements of $X$ the {\em nonzero} letters  in
$Y$.
Define
$ \Bx(n) = \{(a_1,\ldots ,a_n) : a_i \in Y \mbox{ for all }i\}$, the set
of all $n$-tuples of elements of $Y$ (we use $L_0$ instead of $0$ 
for the zero letter for later convenience. We do
not want to confuse the letter $0$ with the vector $0$).
Given $a = (a_1,\ldots ,a_n)\in \Bx(n)$, define the {\em support} of
$a$ by $S(a) = \{i\in [n] : a_i \not= L_0 \}$. 
For $0\leq i \leq
n$, $\Bx(n)_i$ denotes the set of all elements $a\in \Bx(n)$ with
$|S(a)|=i$.
We have 
$$|\Bx(n)|=(|X|+1)^n,\;\;|\Bx(n)_i|=\binom{n}{i} |X|^i.$$
(We take the binomial coefficient $\binom{n}{k}$ to be 0 if $n<0$ or $k<0$).

There is a natural action of the wreath product $G_n$ on $\Bx(n)$ and
$\Bx(n)_i$:
permute the $n$ coordinates followed by
independently acting on the nonzero letters by elements of $G$. In detail,
given $(g_1,g_2,\ldots ,g_n,\pi)\in G_n$ and $a=(a_1,\ldots
,a_n)\in \Bx(n)$, we have
$(g_1,\ldots ,g_n,\pi) (a_1,\ldots ,a_n) =
(b_1,\ldots ,b_n)$, where $b_i = g_i a_{\pi^{-1}(i)}$, if $a_{\pi^{-1}(i)}$
is a
nonzero letter and $b_i = L_0$, if $a_{\pi^{-1}(i)}=L_0$. 
We have the
following direct sum decomposition into $G_n$-submodules of the permutation
representation of $G_n$ on $V(\Bx(n))$:
\beq 
V(\Bx(n)) &=& V(\Bx(n)_0) \oplus V(\Bx(n)_1) \oplus \cdots \oplus
V(\Bx(n)_n).
\eeq
We now introduce some notation. Let $\sigma_1,\ldots ,\sigma_m$,
$\sigma_i\in G^{\wedge}$, be the distinct irreducible $G$-representations
occuring in the multiplicity free $G$-module $V(X)$. We assume that
$\sigma_1$ is the trivial representation. Now enumerate all the elements of
$G^{\wedge}$ as $\sigma_1,\ldots ,\sigma_t$, so that $\sigma_{m+1},\ldots
,\sigma_{t}$ do not appear in $V(X)$. For $i=1,\ldots ,m$, set
$d_i=\dim(V^{\sigma_i})$, so that $d_1=1$ and $d_1+\cdots + d_m = |X|$.

Denote by $\Y_{2,n}(G^{\wedge})$ the set of all $\mu\in \Y_n(G^{\wedge})$
such that

\noi (i) $\mu(\sigma_i)$ is the empty partition, for $i=m+1,\ldots ,t$.

\noi (ii) $\mu(\sigma_i)$ has atmost one part, denoted $p_i(\mu)$, for
$i=2,\ldots ,m$. We have $p_i(\mu)=0$ if $\mu(\sigma_i)$ is the empty
partition. We set $s(\mu)=p_2(\mu)+\cdots +p_m(\mu)$.

\noi (iii) $\mu(\sigma_1)$ has atmost two parts, denoted $a(\mu), b(\mu)$, 
with $a(\mu)\geq b(\mu)$. Just like in item (ii) above, one or both of
$a(\mu),b(\mu)$ may be $0$.

We have the following combinatorial identity (recall that $V^{\mu}$ denotes
the irreducible $G_n$-module parametrized by $\mu\in \Y_n(G^{\wedge})$). 
\bt \label{rti} 
We have
\beq   
(|X|+1)^n 
&=&  \sum_{\mu\in \Y_{2,n}(G^{\wedge})} (1+a(\mu)-b(\mu))\;\dim(V^{\mu}).
\eeq
\et
\pf The proof is in two steps.

\noi (a) Let $C(n,m)$ denote the set of all $m$-tuples of nonnegative
integers with sum $n$. We have, using the multinomial theorem and (\ref{bi})
above,
\beqn 
\lefteqn{(|X|+1)^n}\\
 &=& (d_1+d_2+\cdots +d_m +1)^n \\ 
          &=& (d_2+\cdots +d_m +2)^n \\ 
          &=& \sum_{(p_1,\ldots ,p_m)\in C(n,m)} 
              {\binom{n}{p_1,\ldots ,p_m}}\;d_2^{p_2}\cdots
d_m^{p_m}\;2^{p_1}\\
          &=& \sum_{(p_1,\ldots ,p_m)\in C(n,m)}\; 
              \sum_{k=0}^{\lfloor p_1/2 \rfloor}\; (p_1 -2k+1)\;
              {\binom{n}{p_1,\ldots ,p_m}}\;d_2^{p_2}\cdots
              d_m^{p_m}\;
              \left\{\binom{p_1}{k}-\binom{p_1}{k-1}\right\}.
\eeqn

\noi (b) Let $\mu \in \Y_{2,n}(G^{\wedge})$. We have
\beqn V^{\mu} &=& \oplus_{T\in \tabx(n,\mu)} V_T.
\eeqn 
The dimension of the GZ-subspace $V_T$ of $V^{\mu}$ is clearly
$d_2^{p_2(\mu)}\cdots d_m^{p_m(\mu)}$. 
With $\mu$ bijectively 
associate the pair of 
elements 
$$(a(\mu)+b(\mu),p_2(\mu),\ldots ,p_m(\mu))\in C(n,m) \mbox{ and }
b(\mu)\in \N \mbox{ with }b(\mu)\leq \lfloor (a(\mu)+b(\mu))/2 \rfloor.$$
It is easy to see, using (\ref{dimi}) above, that the
cardinality of $\tabx(n,\mu)$ is  
$${\binom{n}{a(\mu)+b(\mu),p_2(\mu)\ldots
,p_m(\mu)}}\;\left\{\binom{a(\mu)+b(\mu)}{b(\mu)}
-\binom{a(\mu)+b(\mu)}{b(\mu)-1}\right\}.$$

The result now follows from steps (a) and (b) above. $\Box$

We shall now give a representation 
theoretic interpretation to Theorem \ref{rti} above. 

Consider the 
tensor product
$$ \otimes_{i=1}^n V(Y) = V(Y)\otimes \cdots \otimes V(Y)\;\;(n
\mbox{ factors}),$$
with the natural $G_n$-action (permute the factors and then independently
act on the factors by elements of $G$).
There is a $G_n$-linear isomorphism
\beq \label{tp}
V(\Bx(n)) &\cong &
\otimes_{i=1}^n V(Y)
\eeq
given by
$a=(a_1,\ldots ,a_n) \mapsto a_1 \otimes
\cdots \otimes a_n,\;a\in \Bx(n)$.
From now onwards, we shall
not distinguish between $V(\Bx(n))$ and $\otimes_{i=1}^n V(Y)$.
The image of $V(\Bx(n)_i)$ is denoted $(\otimes_{i=1}^n V(Y))_i$.

Consider the canonical decomposition
$$V(X)=W_1\oplus\cdots \oplus W_m,$$
of $V(X)$ into distinct irreducible $G$-submodules, where $W_i$ is
isomorphic to $V^{\sigma_i}$, for $1\leq i \leq m$.
Thus $d_i=\mbox{dim }W_i$, $i=1,\ldots ,m$.

Define the vector $z\in V(Y)$ by
$z=\sum_{x\in X}x$.

For $0\leq i \leq n$
set
\beqn
\Y_{2,n}(G^{\wedge})_i &=& \left\{ \mu \in \Y_{2,n}(G^{\wedge})\;|\;
b(\mu)+s(\mu)\leq i \leq a(\mu)+s(\mu)\right\}.
\eeqn 
\bt \label{rtI}
For $0\leq i \leq n$, $V(\Bx(n)_i)$ is a multiplicity free $G_n$-module
with $G_n$-module isomorphism
\beqn
V(\Bx(n)_i) &\cong & \bigoplus_{\mu\in \Y_{2,n}(G^{\wedge})_i}\; V^{\mu}.
\eeqn
\et
\pf Let $\mu\in \Y_{2,n}(G^{\wedge})$ and $a(\mu)+s(\mu)\leq i \leq
b(\mu)+s(\mu)$. Let $R\in \tabx(n,\mu)$ be as
defined in the proof of Lemma \ref{sl}. We shall exhibit a GZ-subspace $W$ of
$(\otimes_{j=1}^n V(Y))_i$ of type $V_R$, i.e, $W$ is 
closed under the $G^n$-action and, as a $G^n$-module, is isomorphic to
$V^{r_R(1)}\otimes \cdots \otimes V^{r_R(n)}$ and, for $v\in W$ and
$j=1,2,\ldots ,n$, we have
\beq \label{id1}
X_j(v) &=& \frac{|G|}{\dim(V^{r_R(j)})}\;c(b_R(j))v.
\eeq
This will show that $V^{\mu}$ appears in $V(\Bx(n)_i)$. The dimension count
given by Theorem \ref{rti} then completes the proof.

\noi (a)  Set $q=a(\mu)+b(\mu)$. There is an injection
$$\Gamma : V(B(q)) \rar \otimes_{j=1}^q V(Y)$$
given as follows: for $X\seq [q]$, we have $\Gamma(X)=u_1\otimes \cdots
\otimes u_q$, where $u_k=L_0$, if $k\not\in X$ and $u_k=z$, if $k\in X$.

Since $b(\mu)\leq i-s(\mu)\leq a(\mu)$ and $b(\mu) \leq \lfloor
(a(\mu)+b(\mu))/2 \rfloor$, it follows from Theorem \ref{js}
that there is a vector $u\in V(B(q)_{i-s(\mu)})$ (determined uniquely upto
scalars) such that
\beq \label{id2} 
Y_j(u)&=&c(b_R(j))u,\;j=1,\ldots ,q.
\eeq

\noi (b) Let $\sigma \in G^{\wedge}$ and consider the $G_k$-module
$V^{\sigma}\otimes \cdots \otimes V^{\sigma}$ ($k$ factors). It follows
from Theorem \ref{mtl}(i) and (ii)(a) that, for all 
$v\in V^{\sigma}\otimes \cdots
\otimes V^{\sigma}$,
\beq \label{id3}
X_j(v) = (j-1)\frac{|G|}{\dim(V^{\sigma})}v,\;j=1,\ldots ,k.
\eeq

Consider the subspace $W$ of $(\otimes_{j=1}^n V(Y))_i$ given by
$$W= \mbox{Span}(\Gamma(u))\otimes W_2\otimes \cdots \otimes W_2 \otimes
W_3\otimes \cdots \otimes W_3\otimes
\cdots \otimes W_m\otimes \cdots \otimes W_m,$$
where $W_2$ is repeated $p_2(\mu)$ times, $W_3$
is repeated $p_3(\mu)$ times, and so on until $W_m$ is repeated $p_m(\mu)$
times.

Since $g\cdot L_0=L_0$ and $g\cdot z=z$, for all $g\in G$, it follows that
$W$ is closed under the $G^n$-action and, as a $G^n$-module, 
is isomorphic to $V^{r_R(1)}\otimes \cdots \otimes V^{r_R(n)}$.
Moreover, it follows from (\ref{id2}) above that, for $v\in W$,
\beqn
X_j(v) &=& \frac{|G|}{\dim(V^{r_R(j)})}\;c(b_R(j))v,
\;j=1,\ldots ,q.
\eeqn

From (\ref{id3}) above and Theorem \ref{mtl}(v) we see that, for $v\in W$,
\beqn 
X_j(v) &=& \frac{|G|}{\dim(V^{r_R(j)})}\;c(b_R(j))v,\;j=q+1,\ldots ,n.
\eeqn
That completes the proof. $\Box$

{\bf Acknowledgement} The research of the first author was supported by the
Council of Scientific and Industrial Research, Government of India.


\begin{thebibliography}{99.}


\bibitem{cst1} Ceccherini-Silberstein, T., Scarabotti, F., 
Tolli, F., Trees, wreath products, and finite Gelfand pairs,
{\em Adv. in Math} {\bf 206}, 503-537  (2006).

\bibitem{cst2} Ceccherini-Silberstein, T., Scarabotti, F., 
Tolli, F., {\em Representation theory of the symmetric groups, The
Okounkov-Vershik approach, character formulas, and partition algebras }, 
Cambridge University Press (2010).


\bibitem{cst3} Ceccherini-Silberstein, T., Scarabotti, F., 
Tolli, F., {\em Representation theory and harmonic analysis of wreath products of
finite groups}, Cambridge University Press (2014).

\bibitem{f} Filmus, Y.,  
Orthogonal basis for functions over a slice of the Boolean hypercube,
{\em arXiv:} 1406.0142,  (2014).

\bibitem{fkmw} Filmus, Y., Kindler, G., Mossel, E., Wimmer, K., 
Invariance principle on the slice,
{\em arXiv:} 1504.01689,  (2015).


\bibitem{jk} James, G. D., Kerber, A., 
{\em The representation theory of the symmetric groups},
Addison-Wesley (1981).

\bibitem{m} Macdonald, I. G., {\em 
Symmetric functions and Hall polynomials (2nd Edition)},
Oxford University Press (1995).


\bibitem{ms} Mishra, A.,  Srinivasan, M. K., Wreath product action on
generalized Boolean algebras, {\em arXiv:} 1420.8270.



\bibitem{p} Pushkarev, I. A., On the representation theory of wreath
products of finite groups and symmetric groups,
{\em J. Math. Sci.} {\bf 96}, 3590-3599 (1999).


\bibitem{sr1} Srinivasan, M. K., Symmetric chains, Gelfand-Tsetlin chains, 
and the Terwilliger algebra of the binary Hamming scheme,
{\em J. Algebraic Comb.} {\bf 34}, 301-322 (2011).


\bibitem{st} Stanley, R. P.,
{\em Enumerative Combinatorics - Volume 1 (Second Edition)}.
Cambridge University Press (2012).

\bibitem{st1} Stanley, R. P.,
{\em Enumerative Combinatorics - Volume 2}.
Cambridge University Press (1999).




\bibitem{vo1} Vershik, A. M.,  Okounkov, A., A new approach to the
representation theory of symmetric groups. 
{\em Selecta Math. (N.S.) } {\bf 2}, 581-605 (1996).

\bibitem{vo2} Vershik, A. M.,  Okounkov, A., A new approach to the
representation theory of symmetric groups. II. 
{\em J. Math. Sci. (N.Y.)} {\bf 131}, 5471-5494 (2005).

\end{thebibliography}
\end{document}